\IfFileExists{./prepreamble-amspreprint.sty}{\RequirePackage[packages,theorems,changes]{prepreamble-amspreprint}}{}

\documentclass[]{amsart}

\RequirePackage{biblatex}
\RequirePackage{ltxcmds}
\IfFileExists{./preamble-amspreprint.sty}{\RequirePackage[packages,theorems,changes]{preamble-amspreprint}}{}

\usepackage{manuscript}

\IfFileExists{./postpreamble-amspreprint.sty}{\RequirePackage[packages,theorems,changes]{postpreamble-amspreprint}}{}

\makeatletter
\@ifpackageloaded{changes}{
\definechangesauthor[name = {Mateusz Baran}, color = {TolMutedBlue}]{MB}
\definechangesauthor[name = {Ronny Bergmann}, color = {TolMutedRose}]{RB}
}{}
\makeatother

\makeatletter
\@ifpackageloaded{hyperref}{%
	\hypersetup{
		pdftitle = {A modified Riemannian Levenberg-Marquardt Algorithm for robust or constraint optimization on manifolds},
		pdfauthor = {Mateusz Baran, Ronny Bergmann},
		pdfkeywords = {Differential Geometry, Riemannian Geometry, Levenberg-Marquardt, Optimization, Robust Procrustes},
		pdfcreator = {Created using the Scoop Template Engine version 1.6.0.}
	}
}{
	\pdfinfo{
		/Title (A modified Riemannian Levenberg-Marquardt Algorithm for robust or constraint optimization on manifolds)
		/Author (Mateusz Baran, Ronny Bergmann)
		/Subject ()
		/Keywords (Differential Geometry, Riemannian Geometry, Levenberg-Marquardt, Optimization, Robust Procrustes)
		/Creator (Created using the Scoop Template Engine version 1.6.0.)
	}
}
\makeatother

\title[A modified Riemannian Levenberg-Marquardt Algorithm]{A modified Riemannian Levenberg-Marquardt Algorithm for robust or constraint optimization on manifolds}

\author[M. Baran]{Mateusz Baran\orcidlink{0000-0001-9667-5579}}
\address[M. Baran]{AGH University of Krakow 30 Mickiewicz Ave., 30-059 Krakow, Poland}
\email{\detokenize{mbaran@agh.edu.pl}}

\author[R. Bergmann]{Ronny Bergmann\orcidlink{0000-0001-8342-7218}}
\address[R. Bergmann]{Norwegian University of Science and Technology, Department of Mathematical Sciences, NO-7491 Trondheim, Norway and AGH University of Krakow 30 Mickiewicz Ave., 30-059 Krakow, Poland}
\email{\detokenize{ronny.bergmann@ntnu.no}}
\urladdr{https://www.ntnu.edu/employees/ronny.bergmann}

\date{\today}

\dedicatory{}

\begin{document}

\begin{abstract}
We extend the Levenberg-Marquardt method on Riemannian manifolds to a robust variant that
allows to tackle problems from applications where outliers are to be expected.
We formally state the framework for a block-wise variant of the Robust Riemannian Levenberg-Marquardt
algorithm and discuss how known convergence results can be applied here as well.
We further discuss several alternatives for phrasing the sub problem that has to be solved.
Finally we illustrate how the accompanying open source implementation in \manoptjl can be used
to efficiently solve problems such as geodesic regression, Procrustes analysis, subspace Procrustes analysis and bundle adjustment robustly and compare the Levenberg-Marquardt solver
to other solvers for nonsmooth Riemannian optimization.
\end{abstract}


\makeatletter
\ltx@ifpackageloaded{hyperref}{%
\subjclass[2010]{\href{https://mathscinet.ams.org/msc/msc2020.html?t=90C25}{90C25}, \href{https://mathscinet.ams.org/msc/msc2020.html?t=49Q99}{49Q99}, \href{https://mathscinet.ams.org/msc/msc2020.html?t=49M20}{49M20}, \href{https://mathscinet.ams.org/msc/msc2020.html?t=65K10}{65K10}}
}{%
}
\makeatother

\maketitle

%
%
\section{Introduction}%
\label{section:introduction}

Optimization problems on Riemannian manifolds have gained a lot of interest in the last fifteen years.
For smooth optimization, the two books~\cite{AbsilMahonySepulchre:2008,Boumal:2023:1} give a modern
overview on methods like gradient descent, the conjugate gradient method, quasi Newton and trust region based methods.
They are accompanied by the \manopt family~\cite{BoumalMishraAbsilSepulchre:2014:1,TownsendKoepWeichwald:2016:2,Bergmann:2022:1}
of toolboxes to make the algorithms accessible to applications as well.
In the recent years, algorithms designed to solve nonsmooth problems have been generalized to Riemannian manifolds.
Starting with the subgradient method~\cite{FerreiraOliveira:1998:1} and the proximal point method~\cite{FerreiraOliveira:2002:1} as a seminal algorithm, most newer algorithms focus on splitting methods, \ie they exploit the specific structure of the cost function $f\colon \mathcal M \to \mathbb R$.
Prominent examples involving the proximal map are the cyclic proximal point algorithm~\cite{Bacak:2014:1}, the Douglas-Rachford algorithm~\cite{BergmannPerschSteidl:2016:1}, or the Chambolle-Pock algorithm~\cite{BergmannHerzogSilvaLouzeiroTenbrinckVidalNunez:2021:1}.
Based on evaluating the subgradient there are methods like the proximal and convex bundle method~\cite{HoseiniMonjeziNobakhtianPouryayevali:2021:1,BergmannHerzogJasa:2024:1} or the proximal gradient methods~\cite{HuangWei:2021:1,FengHuangSongYingZeng:2021,HuangWei:2023:1,BergmannJasaJohnPfeffer:2025:1,BergmannJasaJohnPfeffer:2025:2}, where some of the approaches have recently also been investigated for possible accelerations~\cite{FengJiangHuangYing:2025}. Another approach where subgradients or proximal maps and gradients are combined are the difference of convex algorithm~\cite{BergmannFerreiraSantosSouza:2024,FerreiraGoncalvesLouzeiroNemethZhu:2026} and the difference of convex proximal point algorithm~\cite{SouzaOliveira:2015,AlmeidaNetoOliveiraSouza:2020}.

Another prominent algorithm to employ a splitting method is the \term{Levenberg-Marquardt algorithm}.
While it was already generalized to Riemannian manifold already by Peeters in 1993~\cite{Peeters:1993}, its theoretical derivation and a convergence analysis was only considered very recently~\cite{AdachiOkunoTakeda:2022}. This algorithm is designed to solve nonlinear least-squares optimization tasks employing the $\ell_2$-norm on a vectorial function.
A related algorithm to solve these problems is the Riemannian Gauß-Newtom algorithm~\cite{BreidingVannieuwenhoven:2018}.
A prominent application is provided by geodesic regression~\cite{FletcherJoshi:2004,Rentmeesters:2011,Fletcher:2013,ShinOh:2022,LoayzaRomeroSibumWelker:2026},
where a best-fitting geodesic in the sense of least squares for time-labelled points on a manifold is computed.
However, in the presence of outliers or certain non-Gaussian noise, \eg Poisson noise, it is usually beneficial to use a more robust norm like the $\ell_1$-norm.
This is for example done in the modelling of nonlinear least squares in the Ceres solver~\cite{CeresSolver:2023}, where they follow the presentation in~\cite{TriggsMcLauchlanHartleyFitzgibbon:2000}, also referred to as “Triggs correction” in~\cite{Zach:2014}.

The authors of Ceres cover the case of optimization on manifolds in the form of “boxplus manifolds”, which are essentially differentiable manifolds with a canonical choice of a retraction.%
\footnote{Beyond Euclidean space these are the special orthogonal groups $\mathrm{SO}(2)$ and $\mathrm{SO}(3)$, the sphere $\mathbb S^n$, and the real projective space $\mathbb{R}\mathbb{P}^n$.}
The effects of curvature of the manifold, like a finite injectivity radius of retractions and selection of Riesz representers for gradients, are largely ignored in their implementation.
They also do not discuss convergence as opposed to~\cite{AdachiOkunoTakeda:2022}.

In this paper we generalise the Riemannian Levenberg-Marquardt algorithm presented in~\cite{AdachiOkunoTakeda:2022} to a robustified and blockwise version.
We also consider operator-based subsolvers for the surrogate model, which avoid explicitly forming Jacobians.
Furthermore, we provide a consistent framework for handling box-constrained problems in Riemannian spaces, a problem relevant in bundle adjustment~\cite{GongMengSeibel:2015}, robot motion optimization~\cite{Prete:2018} and fitting statistical models~\cite{BroughtonCoopeRenaudTappenden:2011,LiMiaoLiWangZhang:2018}.

In this work we consider the \term{Robust Nonlinear Least Squares} problem, which consist of
$m\in \mathbb N$ blocks. The overall problem is of the form
\begin{equation}%
    \label{eq:robust_block_RLM}
    \argmin_{p \in \mathcal{M}} f(p), \qquad\text{ where }
    f(p) = \frac{1}{2}\sum_{i=1}^{m} \rho_i\left(\lVert F_i(p) \rVert_2^2\right),
\end{equation}
where $F_i\colon\mathcal{M} \to \mathbb{R}^{n_i}$ are continuously differentiable functions,
also called \term{residual functions},
and $\rho_i\colon\mathbb{R}_{\geq 0} \to \mathbb{R}_{\geq 0}$ are non-decreasing twice continuously differentiable,
with $\rho_i(0) = 0$ which we refer to as \term{robustifiers}.
\footnote{Triggs et al.~\cite{TriggsMcLauchlanHartleyFitzgibbon:2000} call robustifiers “useful” when they fulfil $\rho'(s) < 1$ as well as $\rho''(s) < 0$ whenever $s$ is in “the outlier region”.}
We further denote the components of $F_i$ by $f_{i,j}$, $j=1,\ldots,n_i$,
\ie $F_i = (f_{i,1}, \ldots, f_{i,n_i})$.
We denote these component functions $f_{i,j} \colon \mathcal{M} \to \mathbb{R}$ mapping into the reals always by small letters.

We further consider box-constrained problems of the form
\begin{equation}%
    \label{eq:robust_block_RLM-box}
    \argmin_{p=(p_{D}, p_{\mathcal{M}}) \in D \times \mathcal{M}} f(p)\qquad\text{ where }
    f(p) = \frac{1}{2}\sum_{i=1}^{m} \rho_i\left(\lVert F_i(p) \rVert_2^2\right),
\end{equation}
where $D$ is the hypercube $[l_1, u_1] \times [l_2, u_2] \times \cdots \times [l_{n_{\mathrm{box}}}, u_{n_{\mathrm{box}}}]$ for some $n_{\mathrm{box}} \geq 0$, where some of the bounds may be infinite and the box part $p_{D}$ is handled using a generalized Cauchy direction method as introduced in~\cite{BaranBergmannPrzybysz:2026}.
As a motivating example we consider a variant of the bundle adjustment problem with bounds constraints on Euclidean parameters such as point positions, camera positions, focal lengths and radial distortion coefficients.

The remainder of this paper is organised as follows: we first introduce the necessary preliminaries and notation in Section~\ref{section:background}.
We derive a novel Riemannian surrogate for the robustified Riemannian Levenberg-Marquardt algorithm in Section~\ref{section:modified_RLM_Surrogate}.
We present the overall robust Riemannian Levenberg-Marquardt algorithm in Section~\ref{section:modified_RLM}, its parameters and especially how to solve the occurring subproblem.
We also discuss its convergence properties.
We extend the framework further to involving bounds constraints in Section~\ref{section:box-constraints}. Numerical examples illustrate properties and advantages of both the model and the implementation in Section~\ref{section:numerical-experiments} and we conclude the paper in Section~\ref{section:conclusion}.

\section{Preliminaries}%
\label{section:background}

In this section, we recall the necessary terms and notions from both Riemannian geometry and optimization. For more details, see for example the textbooks~\cite{DoCarmo:1992:1,AbsilMahonySepulchre:2008,Boumal:2023:1}.
We denote a smooth, finite-dimensional, second-countable, and Hausdorff \term{Riemannian manifold} by $\mathcal M$.
Its dimension is denoted by $d_{\mathcal M}$ or $d$ when the manifold is clear from context.
At every point $p\in\mathcal M$ we denote the \term{tangent space} by $T_p\mathcal M$,
which is a vector space of dimension $d$ and we denote its elements, the \term{tangent vectors}, by $X,Y \in T_p\mathcal M$.
The disjoint union of all tangent spaces yields the \term{tangent bundle}, which we denote by $T\mathcal M$.
The family of inner products $\langle \cdot, \cdot \rangle_p\colon T_p\mathcal M \times T_p\mathcal M \to \mathbb R$ that smoothly vary in $p \in \mathcal M$ is called the \term{Riemannian metric}.
In order to “move away from a point $p$ along a direction $X \in T_p\mathcal M$”
we consider \term{retractions} $R\colon T\mathcal M \to \mathcal M$, which are defined as smooth maps $(p,X) \mapsto R_p(X) \in \mathcal M$ such that each curve $c(t) = R_p(tX)$ satisfies $c(0) = p$ and for its first derivative we have $c'(0) = X$.
A special case of a retraction is the \term{exponential map} $\exp_p \colon T_p\mathcal M \to \mathcal M$, for which the curves $c$ are locally length-minimizing.
The length of these curves also allows to define a distance on the manifold which we denote by $\dM_{\mathcal M}\colon \mathcal M \times \mathcal M \to \mathbb R_{\geq 0}$.
Locally around $0 \in T_p\mathcal M$ the exponential map can be inverted. We call $\log_p(q) \coloneqq \exp_p^{-1}(q)$ the \term{logarithmic map}.
Similarly, one can locally invert a retraction and obtains the \term{inverse retraction} $R_p^{-1}$.
If furthermore each of the curves $c(t) = R_p(tX)$, $p\in\mathcal M$, $X\in T_p\mathcal M$, has zero acceleration at $t=0$,
\ie $c''(0) = 0$, we call $R_p$ a \term{second order retraction}~\cite[Def.~5.42]{Boumal:2023:1}.

For a smooth function $f\colon \mathcal M \to \mathbb R$ we denote its \term{differential}
by $Df\colon T\mathcal M \to \mathbb R$ or point wise as $Df(p)\colon T_p\mathcal M \to \mathbb R$.
We define the \term{Riemannian gradient} $\grad f\colon \mathcal M \to T\mathcal M$ as the Riesz representer of the differential with respect to the Riemannian metric, \ie that for any $p \in \mathcal M$ we have that for all $X \in \mathcal T_p\mathcal M$ it holds
\begin{equation*}
    \Diff{}f(p)[X] = \langle \grad f(p), X\rangle_p.
\end{equation*}
Similarly, we define the \term{Riemannian Hessian} $\Hess f$ as the pointwise linear map $\Hess f(p)\colon T_p\mathcal M \to T_p\mathcal M$ given by the covariant derivative of the gradient
\begin{equation*}
    \Hess f(p)[X] = \nabla_X\grad f,
\end{equation*}
where $\nabla$ denotes the \term{Riemannian connection}.

For a vectorial function $F\colon \mathcal M \to \mathbb R^n$ we can write $F(p) = (f_1(p),\ldots, f_n(p))$ using its component functions $f_j \colon \mathcal M \to \mathbb R$, $j=1,\ldots,n$.
The \term{Jacobian} $\mathcal{J}_F\colon T\mathcal M\to \mathbb R^n$ or \term{differential}
of $F$ is given by
\begin{equation*}
    \mathcal{J}_F(p)[X]
        = \Diff{}F(p)[X]
        = \Bigl( \langle \grad f_j(p), X \rangle_p \Bigr)_{j=1}^{n}.
\end{equation*}
We use calligraphic letters here to emphasize, that the Jacobian is a linear operator independent of a choice of coordinates or basis of the tangent space $T_p\mathcal M$.
Let $p \in \mathcal M$ be fixed. We denote the \term{adjoint} of a linear operator $\mathcal L\colon T_p\mathcal M\to\mathbb R^n$ by $\mathcal L^*\colon\mathbb R^n \to T_p\mathcal M$ and it is defined by
    \begin{equation*}
        \langle \mathcal L^*[x], X \rangle_p
        = \langle x, \mathcal L[X] \rangle,
            \quad \text{ for all } X \in T_p \mathcal{M}, x \in \mathbb{R}^{n},
    \end{equation*}
    where $\langle\cdot,\cdot\rangle$ denotes the default inner product on $\mathbb R^n$.
    Hence the \term{adjoint Jacobian} $\mathcal{J}_F^*\colon\mathbb R$ can be written as
    \begin{equation*}
        \mathcal{J}_F^*(p)[x]
        = \sum_{j=1}^{n} x_j \grad f_j(p),\qquad p \in \mathcal M, x \in \mathbb R^n.
    \end{equation*}
    Similarly, for a fixed point $p \in \mathcal M$ and some $C \in \mathbb R^{n\times n}$ the adjoint of a linear operator $\mathcal L_C[X] \coloneqq C\mathcal J_F(p)[X]$
    is given by by $\mathcal L_C^*[x] = \mathcal J_F^*(p)[C^{\mathrm{T}}x]$.

    For a fixed basis $\{Y_1, \ldots, Y_d\}$ of $T_p \mathcal{M}$ of the tangent space at $p$ we can represent every tangent vector $X \in T_p\mathcal M$
    with unique coefficients $c_l \in \mathbb R$, $l=1,\ldots,d$ as
    \begin{equation*}
        X = \sum_{l=1}^{d} c_l Y_l \in T_p \mathcal{M}.
    \end{equation*}
    Using such a basis $\{Y_1, \ldots, Y_d\}$, we can write $\mathcal{J}_F(p)$ as a matrix
    \begin{equation}
        \label{eq:Jacobian_matrix_in_coordinates}
    J_F(p) =
        \left(
            \langle \grad f_j(p), Y_l \rangle_p
        \right)_{j,l=1}^{n,d}
        \in \mathbb {R}^{n \times d}
    \end{equation}
    and get in short $\mathcal{J}_F(p)[X] = J_F(p) c$, using the coordinates of the tangent vector $X$, \ie the vector $c=(c_1,\ldots,c_d) \in \mathbb R^d$ of coefficients of $X$. We use upper case letters to emphasize that this is a matrix of coefficients representing the Jacobian.

    If the basis $\{Y_1, \ldots, Y_d\}$ of $T_p \mathcal{M}$ is orthonormal, \ie $\langle Y_l, Y_{l'} \rangle_p = 1$ if $l=l'$ and $0$ otherwise, the matrix representation of the adjoint Jacobian is given by the transpose of
    the matrix representation of the Jacobian from~\eqref{eq:Jacobian_matrix_in_coordinates}, that is
    \begin{equation*}
        J_F^*(p) = J_F(p)^{\mathrm{T}} \in \mathbb{R}^{d \times n}.
    \end{equation*}

    Note that depending on the application and the manifold in mind, decomposing a tangent vector with respect to a basis might be expensive to compute. Similar to~\cite[Proposition 1]{AdachiOkunoTakeda:2022}, we therefore prefer to work coordinate-independent in the following.
    For concrete implementations that leaves the freedom to work coordinate-independent or in coordinates.

    Finally, we need the following projection operator.
    Let $n_{\mathrm{box}}\in\mathbb N$ be a natural number and $D = [l_1, u_1] \times [l_2, u_2] \times \cdots \times [l_{n_{\mathrm{box}}}, u_{n_{\mathrm{box}}}]$ be a box in $\mathbb{R}^{n_{\mathrm{box}}}$ and $\mathcal{M}$ a Riemannian manifold. We consider the product manifold with corners $D \times \mathcal{M}$ with the product metric~\cite{Joyce:2010}.
    Let $p = (p_{D}, p_{\mathcal{M}})$ be a point on $D \times \mathcal{M}$.
    The projection $P_{T_{p} (D \times \mathcal{M})} \colon T_p (\mathbb{R}^{n_{\mathrm{box}}} \times \mathcal{M}) \to T_{p} (D \times \mathcal{M})$ is defined as $P_{T_{p} (D \times \mathcal{M})} Y = P_{T_{p} (D \times \mathcal{M})} (Y_{D}, Y_{\mathcal{M}}) = (P_{T_{p_{D}}D} Y_{D}, Y_{\mathcal{M}})$, where
\begin{equation*}
    P_{T_{p_{D}}D} Y_{D,i}
    = \begin{cases}
        Y_{D,i} & \text{ if } l_i < p_{D,i} < u_i \\
        \min\{0, Y_{D,i}\} & \text{ if } p_{D,i} = l_i \\
        \max\{0, Y_{D,i}\} & \text{ if } p_{D,i} = u_i
    \end{cases}
\end{equation*}
for each component $i \in \{1, 2, \dots, n_{\mathrm{box}}\}$.

\section{The Surrogate for Robust Nonlinear Least Squares}%
\label{section:modified_RLM_Surrogate}

In this section we derive a novel linear surrogate for the robust nonlinear least squares problem in block form~\eqref{eq:robust_block_RLM}. Let $p^{(k)} \in \mathcal M$ for some fixed $k\in \mathbb N$ denote a point on the manifold.
We aim to find a damped linear surrogate of the form
\begin{equation}\label{eq:damped_surrogate_model}
s^{(k)}_{\lambda}(X) \coloneqq \frac{1}{2}\lVert \mathcal L_k(X) + y^{(k)}\rVert_2^2 + \lambda \lVert X \rVert_{p^{(k)}}^2,
    \qquad X \in T_{p^{(k)}} \mathcal{M},
\end{equation}
where $\lambda \geq 0$ is the damping factor, and both the linear operator $\mathcal L_k \colon T_{p^{(k)}} \mathcal{M} \to \mathbb{R}^n$ and the vector $y^{(k)} \in \mathbb{R}^n$ have to be determined.
This is inspired by Eq.~(2.1) from~\cite{AdachiOkunoTakeda:2022}.
We write in short $s^{(k)}(X) \coloneqq s_0^{(k)}(X)$ for the undamped surrogate.

We further define
\begin{equation}\label{eq:combined-summand-Gi}
    G_i(p) \coloneqq \frac{1}{2}\rho_i\left(\lVert F_i(p) \rVert_2^2\right),
    \qquad p \in \mathcal{M}, i=1,\ldots,m.
\end{equation}
Then, one way to introduce a surrogate would be to use $G_i$ as $F$ in~\cite{AdachiOkunoTakeda:2022},
and sum these.
In the following we extend this approach and derive a modified surrogate adapting the so-called “Triggs correction”~\cite{TriggsMcLauchlanHartleyFitzgibbon:2000,Zach:2014} to the Riemannian case.

\subsection{The surrogate for a single block}

For this subsection we fix $i \in \{1,\ldots,m\}$ and consider a single summand $F = F_i$ dropping the index
$i$ for $F$, $G$ from~\eqref{eq:combined-summand-Gi}, their range dimension $n=n_i$, and the robustifier $\rho=\rho_i$, respectively.
This also simplifies the component functions, where we write $F(p) = (f_1(p),\ldots,f_n(p)) \in \mathbb{R}^n$.

The goal in the following is to derive the surrogate $s_{i,\lambda}^{(k)}$ of the form~\eqref{eq:damped_surrogate_model} where we also drop the $i$ for the remainder of this subsection.
We introduce the short hand notations $r_k \coloneqq F(p^{(k)})$ for the vector of residuals,
$\mathcal{J}_k \coloneqq \mathcal{J}_F(p^{(k)})$ with adjoint~$\mathcal{J}_k^*$,
as well as $\rho_k \coloneqq \rho(\lVert r_k \rVert_2^2)$,
$\rho_k' \coloneqq \rho'(\lVert r_k \rVert_2^2)$, and
$\rho_k'' \coloneqq \rho''(\lVert r_k\rVert_2^2)$.

For a non-robustified objective function $h(p) \coloneqq \frac{1}{2} \lVert F(p) \rVert_2^2$ the differential reads
\begin{equation*}
    \Diff{}h(p)[X]
    = \langle F(p), \mathcal{J}_F(p)[X] \rangle,
        \quad \text{ for all } X \in T_p \mathcal{M}.
\end{equation*}
Using the adjoint of the Jacobian and a chain rule, we use this to write the Riemannian gradient of $G$ at $p \in \mathcal{M}$ as
\begin{equation}\label{eq:gradG}
    \grad G(p)
    = \rho'\left(\lVert F(p) \rVert_2^2\right)
      \mathcal{J}^*_F(p)\bigl[F(p)\bigr]
    = \rho'\left(\lVert F(p) \rVert_2^2\right)
        \sum_{j=1}^{n} f_j(p) \grad f_j(p)
\end{equation}
Similarly, the Riemannian Hessian of $G$ at $p \in \mathcal{M}$
can be derived using the product rule and the chain rule as
\begin{equation*}
    \begin{split}
        \Hess G(p)[X]
        &= \rho'\left(\lVert F(p) \rVert_2^2\right)
            \Bigl(
            \mathcal{J}^*_F(p)\bigl[ \mathcal{J}_F(p)[X]\bigr]
            + \sum_{j=1}^{n} f_j(p) \Hess f_j(p)[X]
            \Bigr) \\
        &\qquad + 2\rho''\left(\lVert F(p) \rVert_2^2\right)
            \langle F(p), \mathcal{J}_F(p)[X] \rangle
            \mathcal{J}^*_F(p)\bigl[F(p)\bigr].
    \end{split}
\end{equation*}

To derive our surrogate model for $G$, we consider an arbitrary second order retraction $R_p$ and the second order Taylor expansion of $G$
at $p \in \mathcal{M}$ in direction $X \in T_p \mathcal{M}$, cf.~\cite[Sec.~5.9]{Boumal:2023:1}. It reads
\begin{equation}\label{eq:taylor_expansion_G}
    \begin{split}
        G\bigl(R_p(X)\bigr)
        &= G(p)
            + \langle \grad G(p), X \rangle_p
            + \frac{1}{2} \langle \Hess G(p)[X], X \rangle_p
            + o(\lVert X \rVert_p^3)
        \\
        &= G(p)
            + \rho'\left(\lVert F(p) \rVert_2^2\right)
            \langle \mathcal{J}_F^*(p)\bigl[F(p)\bigr], X \rangle_p
        \\
        &\quad + \rho'\left(\lVert F(p) \rVert_2^2\right)
            \Bigl(
                \langle \mathcal{J}_F^*(p)\bigl[ \mathcal{J}_F(p)[X]\bigr], X \rangle_p
                + \frac{1}{2}\sum_{j=1}^{n} f_j(p) \langle \Hess f_j(p)[X], X \rangle_p
            \Bigr) \\
        & \quad + \rho''\left(\lVert F(p) \rVert_2^2\right)
            \langle F(p), \mathcal{J}_F(p)[X] \rangle
            \langle \mathcal{J}_F^*(p)\bigl[F(p)\bigr], X \rangle_p
            + o(\lVert X \rVert_p^3)
    \end{split}
\end{equation}
This extends the work in~\cite{AdachiOkunoTakeda:2022}, where a linear surrogate model was used.
From the second order Taylor expansion, we note that the Hessians $\Hess f_j$ in the second line are usually not available, since they are usually too expensive to compute, both analytically or numerically.
On the other hand, the last term only involves first order information of $F$ ---the Jacobian $\mathcal{J}_F$ and its adjoint--- and where the second derivative of the robustifier $\rho''$ is available and usually cheap to compute.
This motivated the “Triggs correction”, as stated in~\cite[Sec.~4.3]{TriggsMcLauchlanHartleyFitzgibbon:2000} to drop all terms in involving the Hessians $\Hess f_i$ of the residuals but keeping all other summands.

Since we aim to mimic the Taylor expansion of $G$ at $p^{(k)}$ up to second order and the
domain of the Jacobian is $\mathbb R^{n}$ we choose the Ansatz $\mathcal{L}_k(X) = C_k \mathcal{J}_k[X]$ for the linear operator in our surrogate of the form~\eqref{eq:damped_surrogate_model},
where $C_k \in \mathbb{R}^{n\times n}$ is a correction matrix.
The adjoint of this linear operator reads $\mathcal L_k^*(x) = \mathcal J_k^*[C_k^{\mathrm{T}}x]$.

The goal is to determine $y^{(k)} \in \mathbb R^n$ and $C_k$ such that the gradient at the origin $0_p \in T_{p^{(k)}} \mathcal{M}$ of the (undamped) surrogate model $s^{(k)}$ agrees with $G$ and the Hessian of $s^{(k)}$ is the Hessian of $G$ omitting the Hessian terms of the residual components $f_j$, \ie,
\begin{align}
    \grad s^{(k)}(0_p)
    &= \grad G(p^{(k)}) = \rho_k' \mathcal{J}_k^*[r_k],
    \label{eq:grad_surrogate}
    \\
    \label{eq:Hess_surrogate}
    \Hess s^{(k)}(0_p)[X]
    &= \rho_k' \mathcal{J}_k^*[\mathcal{J}_k[X]] + 2\rho_k'' \langle r_k, \mathcal{J}_k[X] \rangle \mathcal{J}_k^*[r_k].
\end{align}

In the following we assume that $r_k$ is not the zero vector. Otherwise we take the contribution to the gradient as the zero vector, and Hessian contribution as $\rho_k' \mathcal{J}_k^*[\mathcal{J}_k[X]]$. If all blocks return zero $r_k$, we finish the algorithm, since we have reached a global minimizer.

We use linearity of the Jacobian and its adjoint to reformulate the Hessian condition as
\begin{equation}\label{eq:Hess_surrogate}
\begin{split}
    \Hess s^{(k)}(0_p)[X]
    &= \mathcal{J}_k^*\Bigl[
        \rho_k' \mathcal{J}_k[X] + 2\rho_k'' \langle r_k, \mathcal{J}_k[X] \rangle r_k
    \Bigr].
    \\
    &=
    \mathcal{J}_k^* \Bigl[
        \Bigl(
            \rho_k' I_n + 2\rho_k'' r_k r_k^{\mathrm{T}}
        \Bigr) \mathcal{J}_k[X]
    \Bigr].
\end{split}
\end{equation}
Note that the term in parentheses is a scaled identity $I_n$ plus a rank-one modification.
From the undamped surrogate model $s^{(k)}$, cf.~\eqref{eq:damped_surrogate_model} together with our Ansatz for~$\mathcal{L}$ we have
\begin{align}
        \grad s^{(k)}(0_p)
        & = \mathcal{J}_k^*[C_k^{\mathrm{T}} y^{(k)}],
        \label{eq:grad_surrogate_form}\\
        \Hess s^{(k)}(0_p)[X]
        & = \mathcal{J}_k^*\bigl[C_k^{\mathrm{T}} C_k \mathcal{J}_k[X]\bigr].
        \label{eq:Hess_surrogate_form}
\end{align}
Combining these two systems of equations~\eqref{eq:grad_surrogate} and~\eqref{eq:Hess_surrogate} with~\eqref{eq:grad_surrogate_form} and~\eqref{eq:Hess_surrogate_form}, we obtain
\begin{equation*}
    C_k^{\mathrm{T}} y^{(k)} = \rho_k' r_k
    \qquad\text{ and }\qquad
    C_k^{\mathrm{T}} C_k
    = \rho_k' I_n + 2\rho_k'' r_k r_k^{\mathrm{T}}.
\end{equation*}
One way to derive a solution for $C_k$ is to take for some
$\alpha_k \in \mathbb{R}$ the Ansatz
\begin{equation}
    \label{eq:C}
    C_k = \sqrt{\rho_k'}(I - \alpha_k P_k),
    \qquad P_k = \frac{r_k r_k^{\mathrm{T}}}{\lVert r_k \rVert_2^2},
\end{equation}%
where we still assume that $r_k$ is not the zero vector.
When we then consider
\begin{equation*}
    C_k^{\mathrm{T}} C_k
    = \rho_k' (I - \alpha_k P_k)^2
    = \rho_k' \Bigl(
        I - 2\alpha_k P_k + \alpha_k^2 P_k
    \Bigr)
    = \rho_k' I
        + \rho_k' \bigl(\alpha_k^2 - 2\alpha_k\bigr) P_k,
\end{equation*}
we obtain a form we have seen before, namely a scaled identity plus a rank-one modification.
Using this form in the two equations~\eqref{eq:Hess_surrogate} and~\eqref{eq:Hess_surrogate_form} we get
\begin{equation*}
C_k^{\mathrm{T}} C_k = \rho_k' I
    + \rho_k' \alpha_k(\alpha_k-2)P_k = \rho_k' I
        + 2\rho_k'' r_k r_k^{\mathrm{T}},
\end{equation*}
which holds if we choose $\alpha_k$ such that
\begin{equation*}
    \rho_k' \alpha_k(\alpha_k-2) = 2 \rho_k'' \lVert r_k \rVert_2^2,
    \Leftrightarrow
    \alpha_k^2 - 2\alpha_k - 2\frac{\rho_k''}{\rho_k'} \lVert r_k \rVert_2^2 = 0.
\end{equation*}
A real solution only exists if $2\rho_k'' \lVert r_k \rVert_2^2 < \rho_k'$. The formula for the smaller of the (possibly) two solutions reads
\begin{equation}
    \label{eq:alpha}
    \alpha_k = 1-\sqrt{1 + 2\frac{\rho_k''}{\rho_k'} \lVert r_k \rVert_2^2}.
\end{equation}
For the root $\alpha_k < 1$, $C_k$ is even positive definite.

To determine $y^{(k)}$, we use the relation from~\eqref{eq:grad_surrogate} and~\eqref{eq:grad_surrogate_form}
that $C_k^{\mathrm{T}} y^{(k)} = \rho_k' r_k$. This linear system is fulfilled by
\begin{equation}
    \label{eq:Y}
    y^{(k)} = \frac{\sqrt{\rho_k'}}{1 - \alpha_k} r_k,
\end{equation}
since $(1-\alpha_k P_k) r_k = (1-\alpha_k) r_k$.

Combining these results, in order to find a minimizer of $s^{(k)}$ from~\eqref{eq:damped_surrogate_model} we obtain the following linear system to solve
\begin{equation}
    \label{eq:surrogate_model_optimality_condition}
    \mathcal{J}_k^* \bigl[ C_k^{\mathrm{T}} C_k \mathcal{J}_k[X] \bigr] = - \mathcal{J}_k^*  [ C_k^{\mathrm{T}} y^{(k)} ],
        \quad \text{ for some } X \in T_{p^{(k)}}\mathcal M,
\end{equation}
also called the normal equations of the surrogate model,
cf.~\cite[Prop.~1]{AdachiOkunoTakeda:2022} or~\cite[Sec.~10.2]{NocedalWright:2006}
for the Euclidean case.
In order for the solution of~\eqref{eq:surrogate_model_optimality_condition}
to both have a unique solution and this solution being a minimizer, we need that
the linear operator on the left hand side is symmetric positive definite.
However, for the existence of the solution for the damped surrogate model $s^{(k)}_\lambda$ with $\lambda > 0$, positive semi-definiteness of the left-hand side matrix is sufficient and can be easily verified to hold.
Following the same motivation as~\cite{Zach:2014}, we also obtain here, that
since
\begin{equation*}
    C_k^{\mathrm{T}}C_k r_k = \rho'_k r_k + 2\rho''_k r_k r_k^{\mathrm{T}} r_k = (\rho'_k + 2\rho''_k\lVert r_k\rVert_2^2)r_k
\end{equation*}
we get that $r_k$ is an eigenvector with eigenvalue $\rho'_k + 2\rho''_k\lVert r_k\rVert_2^2$,
since we assume $r_k$ is not the zero vector.
For any $V$ orthogonal to $r_k$ we get due to $P_kV = \frac{r_k r_k^{\mathrm{T}} V}{\lVert r_k \rVert_2^2} = 0$
that
\begin{equation*}
     \rho'_k(I-\alpha_k P_k)^2V = \rho'_k V
\end{equation*}
hence $\rho_k'$ is an eigenvalue to any eigenvector of the $n-1$ dimensional subspace orthogonal to $r_k$
and hence of multiplicity $n-1$.For numerical stability we put an upper limit on $\alpha$ such that $\alpha < 1-\varepsilon$ for some small $\varepsilon > 0$.
The developers of Ceres further claim that for $\rho''_k < 0$ the choice of $\alpha$ performs poorly
and in that case they the set $\rho''_k=0$, which yields $\alpha = 0$.

Finally, since $s^{(k)}_\lambda$ is defined on a vector space, we can derive its gradient and Hessian at $X \in T_{p^{(k)}}\mathcal{M}$ from classical Euclidean calculus as
\begin{subequations}
    \begin{align}
        \label{eq:grad_damped_surrogate}
        \operatorname{grad} s^{(k)}_\lambda(X) &= \mathcal L_k^*(\mathcal L_k(X) + y^{(k)}) + \lambda X\\
        \label{eq:Hess_damped_surrogate}
        \operatorname{Hess} s^{(k)}_\lambda(X)[Y] &= \mathcal L_k^*(\mathcal L_k(Y)) + \lambda Y,\qquad Y \in T_{p^{(k)}}\mathcal{M}.
    \end{align}
\end{subequations}
Note that also the damped surrogate fulfils the gradient condition~\eqref{eq:grad_surrogate}, since the second term $\lambda X$ vanished at the origin $0_p$,
but we obtain a regularized linear operator form the Hessian, e.g.~when generalising Eq.~\eqref{eq:surrogate_model_optimality_condition} to the damped surrogate.

\subsection{Properties of the overall surrogate}

In this section we combine the result of the each single block from the previous subsection to a surrogate for the overall model in~\eqref{eq:robust_block_RLM}.
For each summand $i \in \{1,\ldots,m\}$ we introduce the notations
$r_{i,k} \coloneqq F_i(p^{(k)}) $,
$\mathcal{J}_{i,k} \coloneqq \mathcal{J}_{F_i}(p^{(k)})$, $\mathcal{J}_{i,k}^*$ its adjoint,
as well as $\rho_{i,k} \coloneqq \rho_i(\lVert r_{i,k} \rVert_2^2)$, $\rho_{i,k}' \coloneqq \rho_i'(\lVert r_{i,k} \rVert_2^2)$,
$\rho_{i,k}'' \coloneqq \rho_i''(\lVert r_{i,k} \rVert_2^2)$, respectively.

We further obtain an $\alpha_{i,k}$ from~\eqref{eq:alpha} for each block and hence obtain both
\begin{equation*}
    C_{i,k} = \sqrt{\rho_{i,k}'}\Bigl(I - \alpha_{i,k} P_{i,k}\Bigr),
    \qquad
    P_{i,k} = \begin{cases}
        \frac{r_{i,k} r_{i,k}^{\mathrm{T}}}{\lVert r_{i,k} \rVert_2^2} & \text{ if } r_{i,k} \neq 0,\\
        0 & \text{ otherwise,}
    \end{cases}
\end{equation*}
and
\begin{equation*}
    y_i^{(k)} = \frac{\sqrt{\rho_{i,k}'}}{1 - \alpha_{i,k}} r_{i,k}.
\end{equation*}
These allow us to write the linear operators
\begin{equation*}
    \mathcal L_{i,k}(X) = C_{i,k}\mathcal J_{i,k}[X]
\quad
\text{ and its adjoint }
\quad
    \mathcal L^*_{i,k}(x) =\mathcal J^*_{i,k}[ C_{i,k}^{\mathrm{T}}x]
\end{equation*}
to build the \term{damped surrogate model for the robustified nonlinear least squares problem}~\eqref{eq:robust_block_RLM} as
\begin{equation*}
    s_{\lambda}^{(k)}(X)
    = \lambda \lVert X \rVert_{p^{(k)}}^2
        + \frac{1}{2}\sum_{i=1}^m\ \bigl\lVert \mathcal L_{i,k}(X) + y_i^{(k)} \bigr\rVert_2^2,
        \qquad \lambda > 0.
\end{equation*}
With $n = \sum_{i=1}^m n_i$
we further introduce the “stacked” operator $\mathcal L_k\colon T_{p^{(k)}}\mathcal M \to \mathbb R^n$, $\mathcal L_k(X) = \bigl(\mathcal L_{1,k}(X),\ldots,\mathcal L_{m,k}(X)\bigr) \in \mathbb R^n$ and the corresponding $y^{(k)} = \bigl(y_1^{(k)},\ldots, y_m^{(k)}\bigr) \in \mathbb R^n$.
The corresponding adjoint $\mathcal L^*_k$ we can write using $x = (x_1,\ldots, x_m) \in \mathbb R^n$, $x_i \in \mathbb R^{n_i}$, as $\mathcal L^*_k(x) = \displaystyle\sum_{i=1}^m \mathcal L^*_{i,k}(x_i)$.

We obtain a form of our surrogate that has the same form as~\cite[Eq.~(2.1)]{AdachiOkunoTakeda:2022} and we further obtain the analogue op their Eq.~(2.2) for our surrogate as
\begin{subequations}
    \begin{align}
    \label{eq:full-surrogate}
    s_{\lambda}^{(k)}(X)
    &= \frac{1}{2}\lVert y^{(k)} + \mathcal L_k(X) \rVert_2^2 + \lambda\lVert X \rVert_{p^{(k)}}
    \\
    \label{eq:split-surrogate}
    &=
    \frac{1}{2}\lVert y^{(k)} \rVert_2^2
    + \langle \mathcal L^*_{k}(y^{(k)}), X \rangle_{p^{(k)}}
    + \frac{1}{2} \lVert \mathcal L_k(X) \rVert_2^2 + \lambda \lVert X \rVert_{p^{(k)}}.
    \end{align}
\end{subequations}
Furthermore, in the second term in Eq.~\eqref{eq:split-surrogate} we obtain combining~\eqref{eq:gradG} and~\eqref{eq:grad_surrogate} that
\begin{equation*}
    \mathcal L^*_{i,k}(y^{(k)}_i)
    = \mathcal J^*[\rho'_{i,k}r_{i,k}] = \rho'_{i,k}\grad F_i(p^{(k)}) = \grad G_i(p^{(k)}),
\end{equation*}
so that overall we have $\mathcal L^*_k(y^{(k)}) = \grad f (p^{(k)})$.

Both sides of the optimality condition~\eqref{eq:surrogate_model_optimality_condition}
decouple into the different blocks $i=1,\ldots,m$, \ie, we obtain
\begin{equation}
    \label{eq:final_optimality_condition}
    \sum_{i=1}^{m} \mathcal{J}_{i,k}^* \bigl[ C_{i,k}^{\mathrm{T}} C_{i,k} \mathcal{J}_{i,k}[X] \bigr]
    = - \sum_{i=1}^{m} \mathcal{J}_{i,k}^*  [ C_{i,k}^{\mathrm{T}} y_i^{(k)} ],
        \quad \text{ for } X \in T_{p^{(k)}}\mathcal M.
\end{equation}
Which for the damped surrogate yields the analogous statement of~\cite[Prop.~1]{AdachiOkunoTakeda:2022} as
\begin{equation}\label{eq:damped-linear-system-short}
    \mathcal L^*_k(\mathcal L_k[X]) + \lambda X
    =
    ( \mathcal L_k \circ \mathcal L^*_k + \lambda \mathcal I_k)(X)
    = - \mathcal L^*_k(y^{(k)})
    = - \grad f(p^{(k)})
\end{equation}
where $\mathcal I_k$ is the identity linear map on $T_{p^{(k)}}\mathcal M$.

For a basis of the tangent space $T_{p^{(k)}} \mathcal{M}$
we can again write this as a linear system. Let
$Y_1,\ldots,Y_d$ be an orthonormal basis of $T_{p^{(k)}} \mathcal{M}$
and $c = (c_1, \ldots, c_d)^{\mathrm{T}} \in \mathbb{R}^d$ the coefficients of some $X \in T_p\mathcal M$.
Then we can write each Jacobian $\mathcal{J}_{i,k}$ as a matrix
$J_{i,k} \in \mathbb{R}^{n_i \times d}$
as in~\eqref{eq:Jacobian_matrix_in_coordinates} and obtain the linear system
\begin{equation}
    \label{eq:final_linear_system}
    \left(
        \sum_{i=1}^{m} J_{i,k}^{\mathrm{T}} C_{i,k}^{\mathrm{T}} C_{i,k} J_{i,k}
    \right) c
    = - \sum_{i=1}^{m} J_{i,k}^{\mathrm{T}} C_{i,k}^{\mathrm{T}} y_i^{(k)}.
\end{equation}
This can also be written in a single matrix $L_kc = z^{(k)}$, where $L_k \in \mathbb R^{d \times d}$ and $z^{(k)} \in \mathbb R^d$ are the two sums from the formula.
Similarly, the damped version just adds a scaled identity matrix $\lambda I_d$ on the left hand side.
Note that while the dimensions of this linear system do not change, the basis $Y_1,\ldots,Y_d$ has to be chosen anew in every iteration of the algorithm,
where the iterate $p^{(k)}$ changes and hence the tangent space $T_{p^{(k)}}\mathcal M$ the linear system~\eqref{eq:final_optimality_condition} is defined on does.

\section{The Modified Riemannian Levenberg-Marquardt Algorithm}%
\label{section:modified_RLM}
In this section we adapt the Riemannian Levenberg-Marquardt algorithm presented in~\cite{AdachiOkunoTakeda:2022} to solve the more general problem~\eqref{eq:robust_block_RLM}.
We describe the algorithm in detail in the following. It is also summarized in Algorithm~\ref{alg:RLM}.

\begin{algorithm}[tbp]
    \caption{A Robust Riemannian Levenberg-Marquardt Algorithm}
    \label{alg:RLM}
    \begin{algorithmic}[1]
        \Require Cost $f$ with components $F_i$, Jacobians $\mathcal J_{F_i}$, and robustifiers $\rho_i$,
        \\
        parameters $\mu_{\mathrm{u}} >\mu_0 \geq \mu_{\mathrm{l}} > 0$, $\beta_{\mathrm{i}} > 1 > \beta_{\mathrm{d}} > 0$, $\eta_{\mathrm{u}} > \eta_{\mathrm{l}} > 0$, $\eta > 0$,
        \\a manifold $\mathcal M$, and a retraction $R$, and a start point $p^{(0)}$.
        \State $k \gets 0$
        \State Initialize $r_{i,k}$, $\rho_{i,k}$, $\rho_{i,k}'$,$\rho_{i,k}''$ for the damped surrogate $s_{\lambda}^{(0)}(X)$
            \Comment{Sec.~\ref{section:modified_RLM_Surrogate}}
        \While{the stopping criterion is not met}
            \State Update $\lambda_k \gets \mu_k\lVert r_k\rVert_2^2$
            \State Compute the minimizer $X^{(k)} \gets \displaystyle\argmin_{X \in T_{p^{(k)}}\mathcal M} s_{\lambda_k}^{(k)}(X)$
            \Comment{Sec.~\ref{subsection:Solving-the-Surrogate}}%
            \label{line:sub_problem}%
            \\[.25\baselineskip]
            \State Compute~$m_{k} \gets \displaystyle\frac{f(p^{(k)}) - f\bigl(R_{p^{(k)}}(X^{(k)})\bigr)}{\frac{1}{2}\bigl(s_{\lambda_k}^{(k)}(0_{p^{(k)}}) - s_{\lambda_k}^{(k)}(X^{(k)})\bigr)}$%
            \\[.25\baselineskip]
            \If{$m_k \geq \eta_{\mathrm{u}}$}\Comment{The model is accurate $\Rightarrow$ allow larger steps}
                \State $\mu_{k+1} \gets
                \max\{\beta_{\mathrm{d}}\mu_k, \mu_{\mathrm{l}}\}$
            \ElsIf{$m_k < \eta_{\mathrm{l}}$}\Comment{The model is inaccurate $\Rightarrow$ reduce step size}
                \State $\mu_{k+1} \gets \min\{\beta_{\mathrm{i}}\mu_k, \mu_{\mathrm{u}}\}$
            \Else\Comment{Otherwise, keep damping parameter}
                \State $\mu_{k+1} \gets \mu_k$
            \EndIf
            \If{$m_k \geq \eta$}\Comment{The model is accurate enough to accept the step}
                \State $p^{(k+1)} \gets R_{p^{(k)}}(X^{(k)})$
                \State Compute $r_{i,k}$, $\rho_{i,k}$, $\rho_{i,k}'$,$\rho_{i,k}''$ for the damped surrogate $s_{\lambda}^{(k)}(X)$
            \Comment{Sec.~\ref{section:modified_RLM_Surrogate}}
            \Else
                \State $p_{k+1} \gets p_k$
                \Comment{Otherwise, keep iterate}
            \EndIf
            \State $k \gets k+1$
        \EndWhile
        \State \Return $p^{(k)}$
    \end{algorithmic}
\end{algorithm}

Besides the ingredients of the problem~\eqref{eq:robust_block_RLM} to solve, we require
an initial model parameter $\mu_0$ and its upper and lower bounds $\mu_{\mathrm{u}}$ and $\mu_{\mathrm{l}}>0$, respectively,
We further use three update parameters: $\eta_{\mathrm{u}} > \eta_{\mathrm{l}} > 0$, and $\eta > 0$, which is slightly more general than stated in~\cite{AdachiOkunoTakeda:2022} and follows the strategy proposed in~\cite{Fan:2006}.
After solving the surrogate subproblem (line~\ref{line:sub_problem}), we evaluate how well this solves the actual problem computing $m_k$.
If the improvement $m_k$ exceeds $\eta_{\mathrm{u}}$, we reduce the damping parameter by multiplying it with $0<\beta_{\mathrm{d}} < 1$, but not below a minimum value of $\mu_{\mathrm{l}}$.
If the improvement $m_k$ is below $\eta_{\mathrm{l}}$, we increase the damping parameter by multiplying it with $\beta_{\mathrm{i}} > 1$, but not beyond a maximum value of $\mu_{\mathrm{u}}$.
We use a separate threshold $\eta$ for accepting a step.

This procedure includes both strategies used in~\cite{AdachiOkunoTakeda:2022} as special cases:
Let their parameters $\beta$ and $\mu_{\min}$ be given and their indicator $\mathrm{flag}^{\mathrm{nz}}$ to indicate whether they expect a nonzero residual or not.
For obtaining their scenario we first set $\eta_{\mathrm{l}} = \eta$, $\mu_{\mathrm{l}} = \mu_{\mathrm{min}}$, $\mu_{\mathrm{u}} = \infty$, and $\beta_{\mathrm{i}} = \beta$.
If their indicator is true, \ie they expect a nonzero residual, we obtain their setup with $\eta_{\mathrm{u}} = \infty$ and $\beta_{\mathrm{d}} = 1$.
The scenario where they do not expect a nonzero residual by setting their $\mathrm{flag}^{\mathrm{nz}}$ to false
is represented here by setting $\eta_{\mathrm{u}} = \eta$ and $\beta_{\mathrm{d}} = \beta^{-1}$.

\subsection{Solving the subproblem}
\label{subsection:Solving-the-Surrogate}

In order to minimize the surrogate model~$s_{\lambda}^{(k)}(X)$ from~\eqref{eq:damped_surrogate_model} in line~\ref{line:sub_problem} of Algorithm~\ref{alg:RLM}
together with the derived results for $y^{(k)}$ and $C_k$ in~\eqref{eq:Y} and~\eqref{eq:C}, respectively,
we have different choices for a subsolver.
In general, these solvers fall into two categories: either using linear operators and an iterative solver
or using a matrix representation of the operators in a basis of the tangent space $T_{p^{(k)}}\mathcal M$.
While the presentation in~\cite{AdachiOkunoTakeda:2022} and their
Matlab code%
\footnote{The authors kindly provided their source code in personal communication.}
indicate, that the authors
were already aware of both approaches, we would like to provide details and advantages of either method.

\subsubsection*{Using linear operators and vector fields: A conjugate residual method}

The right hand side of~\eqref{eq:surrogate_model_optimality_condition} can be interpreted
as a symmetric linear operator on the tangent space, the right hand side as a vector.
Due to our assumption that the residual vector is not the zero-vector, neither is the right hand side.
For the symmetric operator we move to the damped surrogate $s_\lambda^{(k)}$, $\lambda > 0$,
to obtain invertibility. Hence we can employ for example a conjugate gradient residual method~\cite[Alg.~3]{LaiYoshise:2024}%
\footnote{available in \texttt{Manopt.jl}, see~\href{https://manoptjl.org/stable/solvers/conjugate_residual/}{manoptjl.org/stable/solvers/conjugate\_residual/}},
but other iterative solvers that only require evaluations of the linear operator could also be used.
In this approach one does not need to choose a basis of the tangent space.
As an advantage it is neither necessary to assemble the full matrix in coefficients of the basis, especially representing
the Jacobians in a basis like in Eq.~\eqref{eq:final_linear_system}, nor is a reconstruction of the tangent vector
from its coefficients $c$ necessary.

\subsubsection*{Using coordinates}
When choosing a fixed basis for all blocks $F_i$, $i=1,\ldots,m$, allows to
directly operate in this single bases' coordinates.
Having set up the matrix on the right hand side and the vector on the left hand side of Eq.~\eqref{eq:final_linear_system},
one can use arbitrary linear solvers, for example also those exploiting sparsity.
Since any basis of the tangent space can be used, a basis that promotes sparsity or acts as a preconditioner for the linear system can be selected.
Note that such a basis has to be determined anew in every iteration step, since the tangent space $T_{p^{(k)}}\mathcal M$ depends on the current iterate $p^{(k)}$.
In general one can also choose a different basis for every block $F_i$, $i=1,\ldots,m$, to represent
the vector in coordinates. This, however would cause the necessity for a change of basis before
combining them as in Eq.~\eqref{eq:final_linear_system}

\subsection{Convergence}
While our model uses several summands $G_i$, $i=1,\ldots,m$ as introduced in~\eqref{eq:combined-summand-Gi},
the local convergence analysis of~\cite[Sect.~4]{AdachiOkunoTakeda:2022} can be directly adapted to the setting here.
The main adaption to multiple summands and the robustifiers and our slightly more general damping term update only impose minor changes. We summarize these in the following remark

\begin{remark}\label{rem:ConvAdaptions}
We have that the function $F$ from~\cite{AdachiOkunoTakeda:2022} is in our model the vector $\mathbf{G} \coloneqq (G_1,\ldots,G_m)$ and their surrogate $\theta^k$ is replaced by our novel surrogate $s^{(k)}_\lambda$ which includes second order information. This leads to their $F(x_k)$ and their $J_k$ being replaced by our $y^{(k)}$ and $\mathcal L_k$, respectively. We further use a slightly more general rule to update the damping parameter $\lambda_k$ that includes their two scenarios as special cases.
\end{remark}

For completeness we still repeat the main properties and point out differences when they appear in the proofs.

\begin{assumption}[{\cite[Ass.~1]{AdachiOkunoTakeda:2022}}]
    \label{ass:bounded_Jacobian}
    The Jacobian matrix $\mathcal J_F = (\mathcal J_{F_1},\ldots,\mathcal J_{F_m})$
    that maps from the tangent bundle $T\mathcal M$ into $\mathbb R^n$ and its adjoint $\mathcal J^*$ are bounded on the sublevel set $\mathcal N(p^{(0)}) = \{p \in \mathcal M\ |\ f(p) \leq f(p^{(0)})\}$, \ie there exists $M > 0$ such that $\max\{ \lVert \mathcal J(p) \rVert, \lVert \mathcal J^*(p) \rVert \} < M$ holds for all $p \in \mathcal N(p^{(0)})$.
\end{assumption}

This assumption ensures directly that for $\alpha_k < 1$ both  $\mathcal L_k$ and its adjoint are bounded $\mathcal L^*_k$ are bounded by some $\tilde M > 0$.
Note also that because $\rho$ is monotonically increasing, the arguments of $\rho$, $\rho'$ and $\rho''$ are also bounded for successful steps, and thus from twice continuous differentiability of $\rho$ the values of $\rho$, $\rho'$ and $\rho''$ are also bounded at each iteration.
\begin{lemma}
    The $k$th search direction $X^{(k)}$ satisfies
    \begin{equation*}
        s_{\lambda_k}^{(k)}(0_{p^{k}}) - s_{\lambda_k}^{(k)}(X^{k})
    \geq \frac{\lVert \grad f(p^{(k)}) \rVert_{p^{(k)}}^2}{\lVert \mathcal L_k \rVert^2 + \lambda_k}.
    \end{equation*}
\end{lemma}
The proof of this Lemma follows exactly the lines of the proof of~\cite[Lemma~1]{AdachiOkunoTakeda:2022} with the adaptions mentioned in Remark~\ref{rem:ConvAdaptions}.

\begin{lemma}
    The solution $X^{(k)}$ from the subproblem satisfies
        \begin{align*}
        \lVert X^{(k)} \rVert_{p^{(k)}} &\leq \frac{1}{\lambda_k}\lVert\grad f(p^{(k)}) \rVert_{p^{(k)}}
        \intertext{and }
        -\langle \grad f(p^{(k)}), X^{(k)} \rangle_{p^{(k)}} &\geq \frac{\lVert \grad f(p^{(k)}) \rVert_{p^{(k)}}^2}{\lVert \mathcal L_k \rVert^2 + \lambda_k}.
        \end{align*}
\end{lemma}
The proof of this Lemma follows exactly the lines of the proof of~\cite[Lemma~2]{AdachiOkunoTakeda:2022} with the adaptions mentioned in Remark~\ref{rem:ConvAdaptions}.

\begin{lemma}
    Under Assumption~\ref{ass:bounded_Jacobian}, if we have that $\displaystyle\operatorname*{lim\,inf}_{j\to\infty} \lVert \mathbf{G}(p^{(k(j))}) \rVert > 0$ then $\displaystyle\operatorname*{lim\,sup}_{j\to\infty}
    \mu_{k(j)}\lVert X^{(k(j))} \rVert_{p^{(k)}} < \infty$.
\end{lemma}
Besides Remark~\ref{rem:ConvAdaptions} the proof also requires $\tilde M$ instead of $M$ but otherwise is again verbatim the proof of~\cite[Lemma~3]{AdachiOkunoTakeda:2022}. This also holds for the convergence result from~\cite[Theorem~3]{AdachiOkunoTakeda:2022} which for us reads as follows.

\begin{theorem}
    Under Assumption~\ref{ass:bounded_Jacobian} we have $\operatorname*{lim\,inf}_{k\to\infty} \rVert \grad f(p^{(k)}) \rVert_{p^{(k)}}= 0$.
\end{theorem}
We also get the extension to the limes exists when the sequence $\{\mu_k\}$ stays bounded as can be directly adopted from~\cite[Corollary~1]{AdachiOkunoTakeda:2022}.

In complete analogy to this, the local convergence rate analysis from~\cite{AdachiOkunoTakeda:2022} (Lemmata 5--18 and Theorems 2--4) hold here as well with the same adaptions from Remark~\ref{rem:ConvAdaptions} in Notation and constants.

\section{Constrained Robust Nonlinear Least Squares}%
\label{section:box-constraints}

Box constraints can be handled by adapting the dampened quadratic surrogate.
At each iterate $p^{(k)}$, $k = 0, 1, \dots$ we compute the initial search direction $X_k = (X_{D,k}, X_{\mathcal{M},k}) \in T_{p^{(k)}} (D \times \mathcal{M})$
the same way as in the unconstrained case.
Next, we solve the subproblem of finding the first minimum of the damped surrogate $s^{(k)}_{\lambda}$ from Equation~\eqref{eq:damped_surrogate_model}
along the piecewise linear function $d_{\mathrm{PL}}(t) \colon [0, \infty) \to T_{p^{(k)}} (D \times \mathcal{M})$, $d_{\mathrm{PL}}(t) \coloneqq (d_{\mathrm{PL,D}}(t), t X_{\mathcal{M}, k})$ where $d_{\mathrm{PL,D}}$ is defined as
\begin{equation}
	d_{\mathrm{PL,D}}(t)_i = \begin{cases}
		l_i - p_{D,i} & \text{ if } p_{D,i,t} < l_i \\
		t d_{D,k,i} & \text{ if } l_i \leq p_{D,i,t} \leq u_i \\
		u_i - p_{D,i} & \text{ if } p_{D,i,t} > u_i
	\end{cases}
\end{equation}
where $p_{D,i,t} = p_{D,i} + t d_{D,k,i}$ for $i\in \{1, 2, \dots, n_{\mathrm{box}}\}$.
Let the first minimizer of $q_k(t) = s^{(k)}_{\lambda}(d_{\mathrm{PL}}(t))$ be attained at $t_{*,k}$.
The generalized Cauchy direction is then defined as $d_{\mathrm{GCD},k} = d_{\mathrm{PL}}(t_{*,k})$.
To find the minimizer we use the generalized Cauchy direction algorithm~\cite{BaranBergmannPrzybysz:2026}, which is a Riemannian adaptation of the generalized Cauchy point algorithm~\cite{ByrdLuNocedalZhu:1995}.
We then consider $d_{\mathrm{GCD},k}$ as the step proposal instead of the original direction.

Let $\mathcal{H}$ be the Hessian of our surrogate $s^{(k)}_{\lambda}$ defined in Eq.~\eqref{eq:Hess_surrogate}.
The GCD algorithm requires computing the following Hessian values: $\langle e_b, \mathcal{H}[e_b] \rangle_{p^{(k)}}$ and $\langle e_b, \mathcal{H}[X] \rangle_{p^{(k)}}$,
where $e_b$ for $b = 1, 2, \dots, n_{\mathrm{box}}$ is the tangent vector at $p^{(k)}$ that has $b$th standard basis vector as the first component and the zero vector as the second component.
The exact computation is subsolver-dependent.

For the purpose of box constrained minimization, the standard Jacobian $\mathcal{J}_F(p)$ needs to be replaced with the projected Jacobian $\mathcal{J}^P_F(p)$, defined as
\begin{equation}
    \mathcal{J}^P_F(p)[X]
            = \Bigl( \langle P_{T_{p} (D \times \mathcal{M})} \grad f_j(p), X \rangle_p \Bigr)_{j=1}^{n},
\end{equation}
where $P_{T_{p} (D \times \mathcal{M})}$ is the projection defined in Section~\ref{section:background}.
Moreover, stopping criteria should use norm of the projected negative gradient $P_{T_{p} (D \times \mathcal{M})}(-\grad f(p))$ instead of norm of the gradient itself.

\section{Numerical Experiments}%
\label{section:numerical-experiments}

In this section we present examples to illustrate how the modified
Riemannian Levenberg-Marquardt (RLM) algorithm performs in practice.
The algorithm has been implemented in \manoptjl~\cite{Bergmann:2022:1},
version 0.6 or later, %
and we employ manifolds that are implemented in~\manifoldsjl~\cite{AxenBaranBergmannRzecki:2023} version 0.11.28.
The experiments were performed using Julia 1.12.6
and especially regarding time measurements on a machine running Linux Mint 22.3 on an AMD Ryzen 9 9950X3D CPU.
\subsection{Geodesic Regression}
\label{sec:RobustGeodesicRegression}
Following~\cite[Sec.~3]{Fletcher:2013}, we consider the task of \term{gedoesic regression}
but extending it here to a robust variant.
Let $\mathcal M$ be a Riemannian manifold and let $q_1,\ldots,q_m \in \mathcal M$, $m \in \mathbb N$, be points
thereon associated to scalar values $t_1,\ldots,t_m \in \mathbb R$.
We further denote the geodesic starting in a point $p \in \mathcal M$ in direction $X \in T_p\mathcal M$ by
\begin{equation*}
    \gamma_{p,X}(t) \coloneqq \exp_p(tX).
\end{equation*}
The task of geodesic regression is given as the problem on the tangent bundle $T\mathcal M$
to find the geodesic $\gamma_{p,X}$ that “best fits” the data, \ie such that $\gamma_{p,X}(t_i) \approx q_i$, $i=1,\ldots,m$,
where we still have the freedom to choose in what sense we mean to address the approximately equal condition.
We consider
\begin{equation*}
    \argmin_{(p,X) \in T\mathcal M} \sum_{i=1}^m \rho_i\bigl( \dM_{\mathcal M}^2(\gamma_{p,X}(t_i), q_i)\bigr),
\end{equation*}
where $\rho_i$, $i=1,\ldots,m$, are robustifiers to be chosen.
If we set $\rho_i(x) = \sqrt{x}$ we obtain the case of \term{robust geodesic regression}, see also~\cite{ShinOh:2022}, where they use a gradient descent, though the function is not differentiable everywhere.

Here we restrict ourselves to use the same robustifier $\rho_i = \rho$ for all summands.
Setting $\rho = \rho_{\mathrm{lsq}}$, $\rho_{\mathrm{lsq}}(x) = x$,
we obtain the classical least squares objective from~\cite{Fletcher:2013}.
This choice fits very well when the data points $q_i$ are for example perturbed by Gaussian noise.

If, however, the data has outliers, it is beneficial to minimize the 1-norm instead.
Hence we use a robustifier that approximates the 1-norm while being differentiable, the so-called \term{Huber function}
\begin{equation*}
 \rho_{\mathrm{Hu}}\colon\mathbb R_{\geq 0} \to \mathbb R_{\geq 0}
\qquad
\rho_{\mathrm{Hu}}(x) \coloneqq \begin{cases} x & \text{ if } x \leq 1\\ 2\sqrt{x}-1 & \text{ if } x > 1.\end{cases}
\end{equation*}
We can scale any robustifier $\rho$ using%
\footnote{This form indeed scales the single functions $f_i$, or here the distances $\dM_{\mathcal M}$, by a factor $a$.}
$s_{\rho, a}(x) \coloneqq a^2 \rho(\frac{x}{a^2})$
and obtain the scaled Huber robustifier $\rho_{\mathrm{Hu},a} = s_{\rho_{\mathrm{Hu}},a}$.
This keeps the quadratic term close to zero, namely for $x \leq a^2$ or in our objective for distances less or equal than $a$, and uses a linear term otherwise.

\begin{figure}[tbp]
    \centering
    \begin{tikzpicture}
  \begin{axis}[
    name=data, width=0.48\textwidth, axis equal image, scale only axis,
    enlargelimits=false, axis lines=none,
    xmin=0, xmax=1, ymin=0, ymax=1,
  ]
    \addplot graphics[xmin=0,xmax=1,ymin=0,ymax=1] {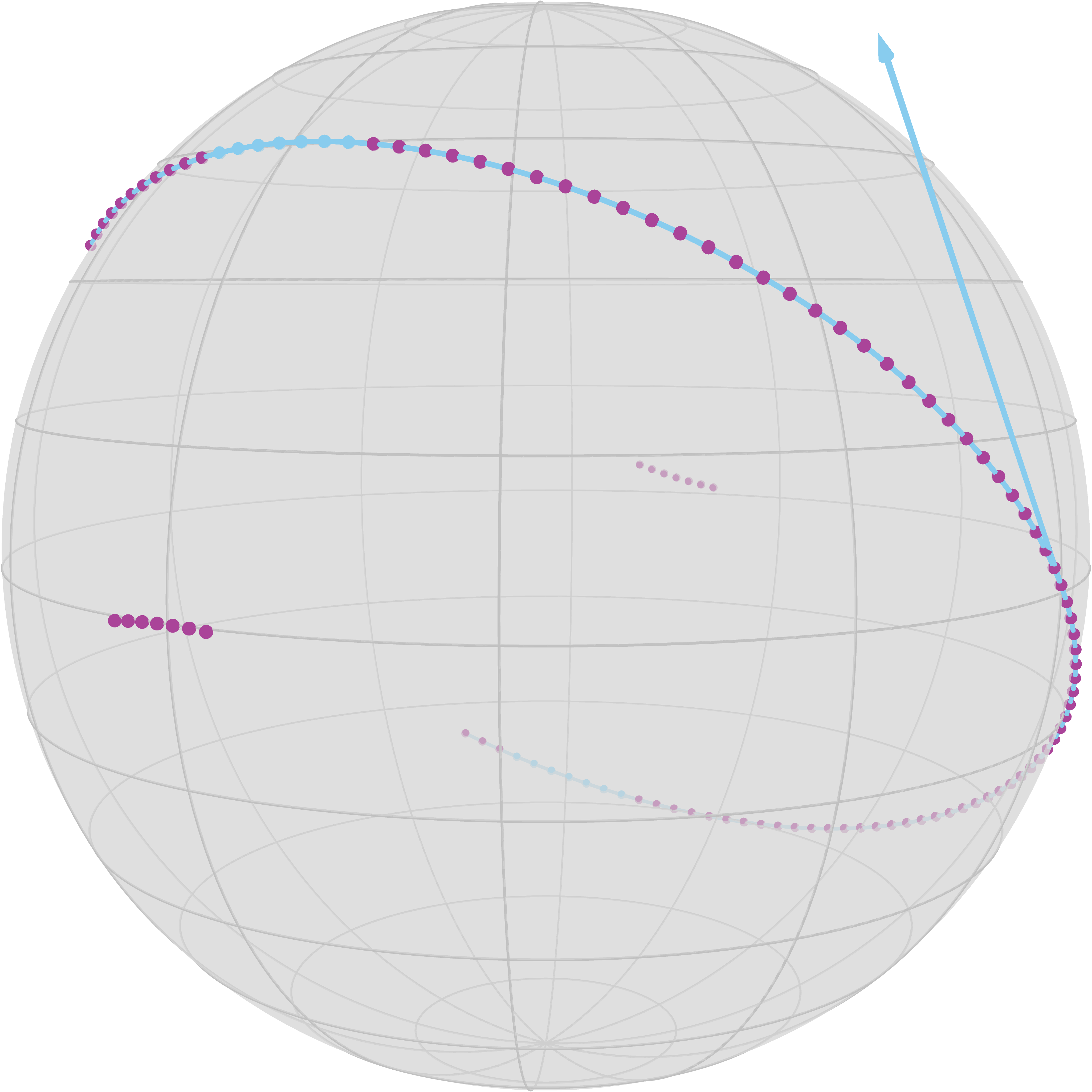};
  \end{axis}

  \begin{axis}[
    name=result,
    at={(data.outer north east)},
    anchor=outer north west,
    xshift=0.02\textwidth,
    width=0.48\textwidth,
    scale only axis,
    axis equal image,
    enlargelimits=false,
    axis lines=none,
    xmin=0, xmax=1,
    ymin=0, ymax=1,
  ]
    \addplot[forget plot] graphics[xmin=0,xmax=1,ymin=0,ymax=1] {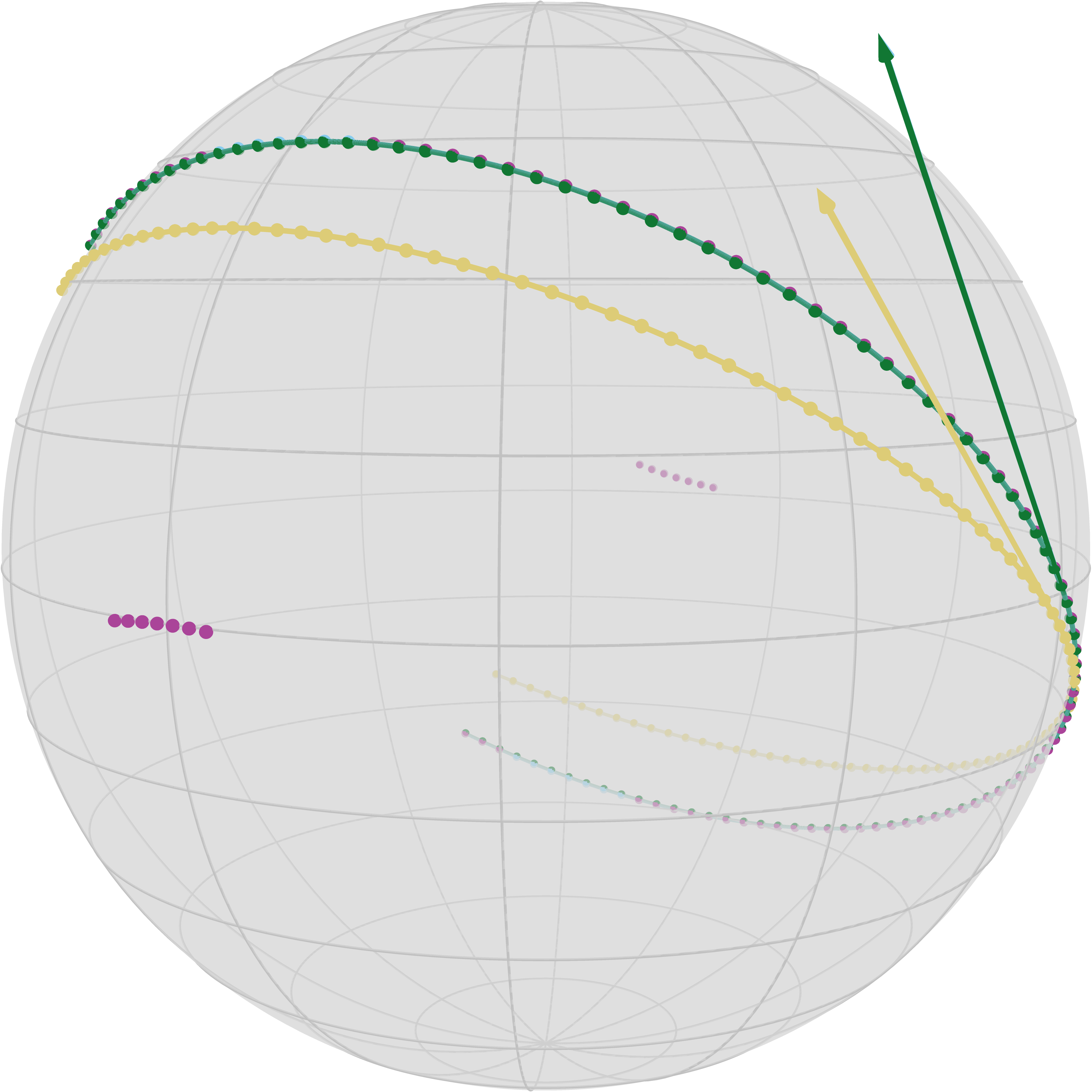};
  \end{axis}
  \node[
    anchor=north,
    at={($(data.outer south)!0.5!(result.outer south)$)},
    yshift=0.25em,
  ] {
    \begin{tikzpicture}
      \begin{axis}[
        hide axis,
        scale only axis,
        width=0pt,
        height=0pt,
        trim axis left,
        trim axis right,
        xmin=0, xmax=1, ymin=0, ymax=1,
        legend style={
          draw=none, font=\small, legend columns=4, column sep=1.125em,
        },
        legend cell align=left,
      ]
    \addlegendimage{TolMutedCyan, very thick, sharp plot, mark=*}
    \addlegendentry{\hspace{-1em}geodesic $\gamma_{p,X}$}
    \addlegendimage{TolMutedPurple, very thick, sharp plot, mark=*, only marks}
    \addlegendentry{\hspace{-1em}Data $q_i$}
    \addlegendimage{TolMutedSand, very thick, sharp plot, mark=*}
    \addlegendentry{\hspace{-1em}Least Squares $\gamma_{p^\ast,X^\ast}$}
    \addlegendimage{TolMutedGreen, very thick, sharp plot, mark=*}
    \addlegendentry{\hspace{-1em}Robust $\gamma_{p^{\star},X^{\star}}$}
      \end{axis}
    \end{tikzpicture}
  };
\end{tikzpicture}
\caption{The original geodesic $\gamma_{p,X}$ (cyan) and data $q_i$, $i=1,\ldots,100$ with outliers (purple).
This yields a reconstruction with least squares regression geodesic $\gamma_{p^\ast,X^\ast}$ (sand, right) that is influenced by
the outliers, while the robust regression geodesic $\gamma_{p^\star,X^\star}$ (green) reconstructs the original geodesics.
For visualization purposes are the tangent vectors scaled by $\frac{1}{2}$.
}%
\label{fig:RobustRegressionSphere}%
\end{figure}

As a concrete example, we take the sphere $\mathcal M = \mathbb S^2$ and generate the data
from the geodesic $\gamma_{p,X}(t)$ with $p = (0,1,0)^{\transp}$ and $X = \frac{\pi}{2}(1,0,1)^{\transp}$,
which we sample at $m=100$ equidistant points for $t$ on $[-1,1]$.
We then generate two sets of outliers for the indices $i \in [4,\ldots,10]\cup[83,\ldots,89]$,
by rotating the tangent vector $Y_i = \log_{q_i}p$ by 90 degrees anticlockwise and
moving the point $q_i$ a length of $r=\frac{\pi}{2}$ into this direction.
The original geodesic and the data with outliers are show in Figure~\ref{fig:RobustRegressionSphere}
on the left.
The first set of outliers is in the back and moved upwards, while the second set in the front
is moved downwards by the length $r$.

We compare two runs of Algorithm~\ref{alg:RLM}, both using a conjugate gradient residual subsolver~\cite[Alg.~3]{LaiYoshise:2024}:
a least squares $\rho = \rho_{\mathrm{lsq}}$ and a robust $\rho = \rho_{\mathrm{Hu},a}$, $a = 10^{-4}$,
model.
Both use the same parameters for the algorithm
\begin{equation*}
    \eta_{\mathrm{u}} = \frac{1}{2},\quad
    \eta_{\mathrm{l}} = \eta = \frac{1}{5},
    \qquad
    \beta_{\mathrm{i}} = 8,\quad \beta_{\mathrm{d}} = \frac{1}{8},
    \qquad
    \mu_0 = \mu_{\mathrm{l}} = 10^{-5},\quad \mu_{\mathrm{u}} = \infty,
\end{equation*}
and the Riemannian center of mass of the data $q_i$ as initial $p^{(0)}$ as well as $X^{(0)} = \log_{p^{(0)}} q_m$ as initial tangent vector.

Key figures of both solver runs are collected in Table~\ref{tab:RobustRegression}.
\begin{table}
    \caption{Key figures from the Least Squares and Robust geodesic Regression on the Sphere}
    \label{tab:RobustRegression}
    \begin{tabular}{rclll}
    \toprule
    & \# Iterations & final cost & mean squared error to $\gamma_{p,X}(t_i)$\\
    \midrule
    Least Squares & 169 & $3.3547$ & $1.3904\times10^{-2}$\\
    Robust & 116 & $1.0995\times10^{-3}$ & $2.2737\times10^{-6}$\\
    \bottomrule
\end{tabular}
\end{table}
The number of iterations might depend also on the actual algorithm parameters chosen, but they do not vary much.
However, both the final cost and the mean distance to the true data (without outliers) is much lower for the robust case.

This can also be seen in Figure~\ref{fig:RobustRegressionSphere} on the right,
where we add both results additionally to the original data and the measured (outlier) data.
The least squares regression geodesic (sand) is skewed to one side by the outliers,
while the robust geodesic (green) is visually indistinguishable both in the sample points
as well as the obtained minimizer $p^{\star}, X^{\star}$.

\subsection{Robust and Robust Subspace Procrustes}
Let two data sets $A \in \mathbb{R}^{d\times n}$,  $A = (a_1,\ldots,a_n)$,
and $B \in \mathbb R^{k\times n}$, $B = (b_1,\ldots,b_n)$,
of measured data be given, where we obtain a set of columns of values $a_i \in \mathbb R^d$, $b_j \in \mathbb R^k$ each, where $k \leq d$. This could be interpreted as $n$ experiments measuring $d$ values, but where for $B$ with $k<d$ the data is in some sense incomplete.

\subsubsection*{The robust case $k=d$ on the rotation matrices $\mathrm{SO}(d)$}
For the case where the number of measurements agrees, \ie $k=d$ the \term{Procrustes Problem}
aims to find the best rotation $p \in \mathrm{SO}(d)$,
$\mathrm{SO}(d) = \{ p \in \mathbb R^{d\times d} \,\mid\, p^{\mathrm{T}}p = I_d, \det(p) = 1\}$,
such that $A \approx pB$, where the classical problem uses the Frobenius norm
to measure what “best” means, \ie,
\begin{equation*}
    \argmin_{p \in \mathrm{SO}(d)}\ \lVert A - pB \rVert_{\mathrm{F}}^2,
    \qquad
    \lVert C \rVert_{\mathrm{F}}^2 \coloneqq \sum_{i,j=1}^{d,n} \lvert c_{i,j} \rvert^2.
\end{equation*}
This problem has a closed form solution based on the singular value decomposition (SVD)
and was introduced in~\cite{Wahba:1965}.
In~\cite{JasaBergmannKuemmerleAthreyaLubberts:2025} the authors looked at a
\term{robust Procrustes problem}
\begin{equation}\label{eq:robustProcrustesCost}
    \argmin_{p \in \mathrm{SO}(d)}\ \lVert A - pB \rVert_{\mathrm{R}},
    \qquad
    \lVert C \rVert_{\mathrm{R}} \coloneqq \sum_{i=1}^{n} \lVert c_i \rVert_2,
\end{equation}
where $c_i$, $i=1,\ldots,n$ again denote the columns of the matrix $C = (c_1,\ldots,c_n) \in \mathbb R^{d\times n}$.
This is inspired by so-called mixed norms often appearing in signal and image processing~\cite{Kowalski:2009,EldarMishali:2009}.
The authors in~\cite{JasaBergmannKuemmerleAthreyaLubberts:2025} used the derivative-free LTMADS algorithm~\cite{Dreisigmeyer:2007} to solve \eqref{eq:robustProcrustesCost}.

In this section we consider the following problem with a single robustifier $\rho$ used in all summands and solve
\begin{equation}\label{eq:robustProcrustesCostHuber}
    \argmin_{p \in \mathcal M} \sum_{i=1}^d\rho \bigl( \lVert a_i - pb_i \rVert_{2}^{2}\bigr)
\end{equation}
for the case $\mathcal M = \mathrm{SO}(d)$.
Within a single summand, the residual function reads $F_i(p) = a_i - pb_i$.
Its differential $\mathcal J_{F_i}(p)\colon \mathfrak{so}(d) \to \mathbb R^n$ is the classical (Euclidean) differential rephrased on the Lie algebra $\mathfrak{so}(d)$, \ie, the set of skew symmetric matrices.
given by
\begin{equation*}
    \mathcal J_{F_i}(p)[X] = DF_i(p)[X] = -pXb_i
    \qquad
    \text{ and }
    \qquad
    \mathcal J^*_{F_i}(p)[y] = -\mathrm{skew}(p^{\mathrm{T}}yb_i^{\mathrm{T}}),
\end{equation*}
where $\mathrm{skew}(A) = \frac{1}{2}(A-A^{\mathrm{T}})$ is the projection onto the set of skew symmetric matrices.
Here the classical Euclidean adjoint is projected, and we again have to take into account, that we represent the result in the Lie algebra.

\pgfplotstableread[col sep=comma]{data/SOd.csv}\SOddatatable
\pgfplotstabletranspose[string type, columns={dim,n}, colnames from=d]\SOdexperiments\SOddatatable
\begin{table}[tbp]
    \pgfplotstabletypeset[
    string type,
    column type={c},
    every head row/.style={before row=\toprule,after row=\midrule},
    every last row/.style={after row=\bottomrule},
    columns/colnames/.style={column name={$d$}, column type={l}},
    every row 0 column 0/.style={string type,
        postproc cell content/.code={%
      \pgfkeyssetvalue{/pgfplots/table/@cell content}{$\mathrm{dim}_{\mathrm{SO}(d)}$}}},
    every row 1 column 0/.style={string type,
        postproc cell content/.code={%
      \pgfkeyssetvalue{/pgfplots/table/@cell content}{$n$}}},
    ]\SOdexperiments
    \caption{Experiments (as columns) for the robust Procrustes problem.}
    \label{table:SOd}
\end{table}

We run the experiment on data for $\mathcal M = \mathrm{SO}(d)$, $d=3,\ldots,15$.
For each $d$ we create a matrix $A$ with $n=\frac{d(d+1)}{2}$ columns with all entries in $[0,1]$.
These parameters are summarized for the thirteen experiments in Table~\ref{table:SOd}.
We then create $B$ as a copy of $A$, where we create four outliers by adding $\frac{1}{10}$ to four entries.
We then generate a random rotation matrix $p^*$
by using \texttt{rand(M; vector\_at = $I_d$, $\sigma$=0.5/d)} to generate a random tangent vector in the tangent space of the identity matrix $I_d$ using the random functions from \manifoldsjl. We apply the (inverse) rotation $(p^*)^{\transp}$ to generate $B$.
This also implies that we have $\lVert A - p^*B\rVert_{\mathrm{R}} = \frac{4}{10}$ for all $d$ by construction. This is at least a local minimizer, though there might exist better local minimizers, since globally the function is nonconvex on $\mathrm{SO}(d)$.

For the robust RLM algorithm, we again use the Huber robustifier $\rho = \rho_{\mathrm{Hu},a}$, here with a scaling of $a = 10^{-5}$. We further set $\beta_{\mathrm{i}} = 4$, $\eta = \frac{1}{5}$, and the initial damping to $\mu_0 = 1\cdot10^{-7}$.
We further set the threshold $\varepsilon = 10^{-4}$ for the scaling $\alpha_k$ in~\eqref{eq:alpha} and use the strict scaling rule from Ceres, i\,e.\ we set $\rho_k'' = 0$ whenever the second derivative of the robustifier $\rho_k''$ is negative.
We compare this to running LTMADS, a mesh adaptive direct search algorithm adapted from~\cite{Dreisigmeyer:2007} available in \manoptjl%
\footnote{see \href{https://manoptjl.org/stable/solvers/mesh_adaptive_direct_search/}{manoptjl.org/stable/solvers/mesh\_adaptive\_direct\_search/}},
that does not require derivative information. It can hence be used to directly minimize~\eqref{eq:robustProcrustesCost}. We use the default parameters from \manoptjl just increasing the maximal number of iterations to $2\cdot10^{4}$. We use the identity matrix, $p^{(0)} = I_d$ as a start point for both algorithms and measure the runtime using \chairmarksjl with 5 samples, 3 evaluations each.
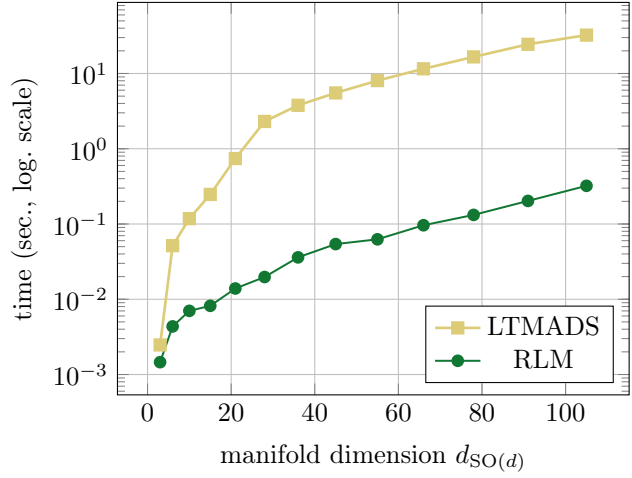
\begin{figure}[tbp]
    \centering
    \begin{tikzpicture}
        \begin{axis}[
            xlabel={manifold dimension $d_{\mathrm{SO}(d)}$}, ylabel={time (sec., log.~scale)},
            legend pos=south east, grid = major,
            ymode = log,
            height=.33\textheight,
            width=.66\textwidth,
        ]
            \addplot[TolMutedSand, line width=1pt, mark size = 2, mark=square*] table [x=dim, y=t2] {\SOddatatable};
            \addlegendentry{LTMADS}
            \addplot[TolMutedGreen, line width=0.75pt, mark size = 2, mark=*] table [x=dim, y=t1] {\SOddatatable};
            \addlegendentry{RLM}
        \end{axis}
    \end{tikzpicture}
    \caption{Runtime comparison of our RLM (green) to the LTMADS algorithm (sand) on different special orthogonal groups.}%
    \label{fig:SOd-runtime}
\end{figure}
The resulting times are shown in Figure~\ref{fig:SOd-runtime}.
While the robust RLM algorithm finishes always in less that $0.322$ seconds, LTMADS needs a lot of iterations, even hits the increased maximal number of iterations starting from $d=9$ and runs for about $32.3$ seconds for the largest experiment $d=15$.
RLM is always faster, for the larger experiments, $d > 7$, at least by a factor $100$.
\begin{figure}[tbp]
    \centering
    \begin{tikzpicture}
    \begin{axis}[
        xlabel={manifold dimension $d_{\mathrm{SO}(d)}$},
        ylabel={final cost $f(p) = \lVert A - pB\rVert_{\mathrm{R}}$},
        legend pos=north west, grid = major,
        height=.33\textheight,
        width=.66\textwidth,
    ]
        \addplot[TolMutedSand, line width=1pt, mark size = 2, mark=square*] table [x=dim, y=f2] {\SOddatatable};
        \addlegendentry{LTMADS}
        \addplot[TolMutedGreen, line width = 0.75pt, mark size = 2, mark=*] table [x=dim, y=f1]{\SOddatatable};
        \addlegendentry{RLM}
        \end{axis}
    \end{tikzpicture}
    \caption{Cost value of the resulting special orthogonal matrix
    of our RLM (green) to the LTMADS algorithm (sand).
    }%
    \label{fig:SOd:Cost}
\end{figure}
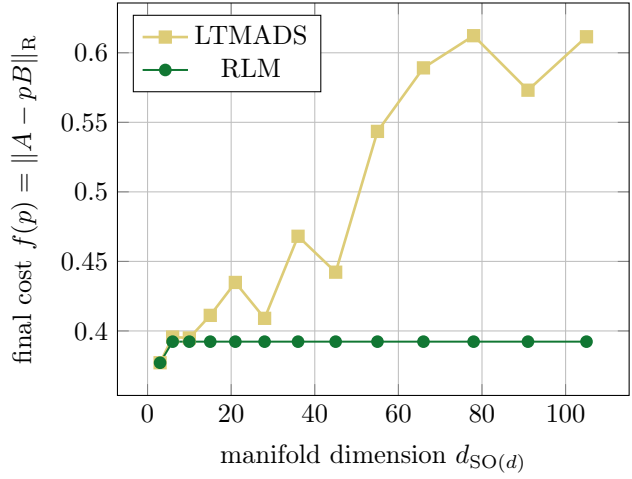

While RLM solves the objective~\eqref{eq:robustProcrustesCostHuber}, \ie, a smoothened version of
the original objective~\eqref{eq:robustProcrustesCost}, LTMADS tackles the
original, nonsmooth objective $f(p) = \lVert A - pB\rVert_{\mathrm{R}}$ directly.
We still compare the final cost value with respect to this nonsmooth objective $f$ of both resulting final iterates in Figure~\ref{fig:SOd:Cost}.
For $d=3$ both find a rotation that provides an even better
alignment than $f(p^*) = 0.4$.
After that, RLM always finds a local minimizer slightly less in cost (a value of about $0.3923$) than $p^*$, while LTMADS slowly deviates
from this value. Starting from $d=9$, where LTMADS hits its maximal number of iterations, the deviation becomes more prominent.

\subsubsection*{The incomplete data and robust case $k<d$ on the Stiefel manifold $\mathrm{St}(d,k)$}
If the measurements of $B$ are incomplete, or in other words $k < d$,
we can still look for the best orthonormal frame $p \in \mathrm{St}(k,d)$, where
$\mathrm{St}(d,k) \coloneqq \{ p \in \mathbb R^{d\times k}\,\mid\, p^{\mathrm{T}}p = I_k\}$
denotes the \term{Stiefel manifold}, \ie the Riemannian manifold of all bases of $k$-dimensional subspaces of $\mathbb R^d$.
This also motivates the phrasing “incomplete” data, since by requiring that $A \approx pB$, we are asking for a “best fitting” orthonormal basis of a $k$-dimensional subspace of $\mathbb R^d$.
When considering classically the Frobenious norm again, this is also called \term{unbalanced Procrustes}~\cite{FulovaTrnovska:2023}.
With the robust norm $\lVert\cdot\rVert_{\mathrm{R}}$ from before, we obtain a \term{robust subspace Procrustes problem}
\begin{equation*}
    \argmin_{p \in \mathrm{St}(d,k)}\ \lVert A - pB \rVert_{\mathrm{R}},
\end{equation*}
which can be interpreted as follows: We try to find an orthonormal basis $p_1,\ldots,p_k \in \mathbb R^d$ that “maps” the columns of $B$ to the columns of $A$.
This is the same objective as in~\eqref{eq:robustProcrustesCostHuber} but now with $\mathcal M = \mathrm{St}(d,k)$.

On the Stiefel manifold, the representation of tangent vectors is done in the embedding, \ie using the Euclidean space $\mathbb R^ {d\times k}$.
The Jacobian is hence just the Euclidean Jacobian
and its adjoint is obtained by projecting the Euclidean adjoint $-yb_i^{\mathrm{T}}$ to the tangent space
\begin{equation*}
    \mathcal J_{F_i}(p)[X] = DF_i(p)[X] = -Xb_i
    \qquad
    \text{ and }
    \qquad
    \mathcal J^*_{F_i}(p)[y] = -yb_i^{\mathrm{T}} + p\cdot\mathrm{sym}(-p^{\mathrm{T}}yb_i^{\mathrm{T}}),
\end{equation*}
where $\operatorname{sym}(A) = \frac{1}{2}(A+A^{\mathrm{T}})$ is the projection onto the set of symmetric matrices.

We run the experiment with the same settings and the same data generation as described in the last section.
We just have to add a final step in the generation of the data $B$ which is reduced in to (its first) $k < d$ rows.
Similarly the optimal value $p^* \in \mathrm{St}(d,k)$ is obtained from the random rotation from the last experiment by removing the last $d-k$ columns.
The parameters with which experiments run are summarized in Table~\ref{table:Stdk}. Compared to the previous experiment we additionally have to choose the reduced data dimension $k$, which is set it to $k = d-\lceil \frac{d}{3} \rceil$ \ie to be less or equal to 2/3rds of the measurements.
Furthermore, on the Stiefel manifold, the exponential map is expensive to compute.
We use the default from \manifoldsjl, namely the polar retraction, since its differential also provides a vector transport as well, which is required to run LTMADS.

\pgfplotstableread[col sep=comma]{data/Stdk.csv}\Stdkdatatable
\pgfplotstabletranspose[string type, columns={k,dim,n}, colnames from=d]\Stdkexperiments\Stdkdatatable
\begin{table}[tbp]
    \pgfplotstabletypeset[
    string type,
    column type={c},
    every head row/.style={before row=\toprule,after row=\midrule},
    every last row/.style={after row=\bottomrule},
    columns/colnames/.style={column name={$d$}, column type={l}},
    every row 0 column 0/.style={string type,
        postproc cell content/.code={%
      \pgfkeyssetvalue{/pgfplots/table/@cell content}{$k$}}},
    every row 1 column 0/.style={string type,
        postproc cell content/.code={%
      \pgfkeyssetvalue{/pgfplots/table/@cell content}{$\mathrm{dim}_{\mathrm{St}(d,k)}$}}},
    every row 2 column 0/.style={string type,
        postproc cell content/.code={%
      \pgfkeyssetvalue{/pgfplots/table/@cell content}{$n$}}},
    ]\Stdkexperiments
    \caption{Experiments (as columns) run for incomplete robust Procrustes problem.}
    \label{table:Stdk}
\end{table}

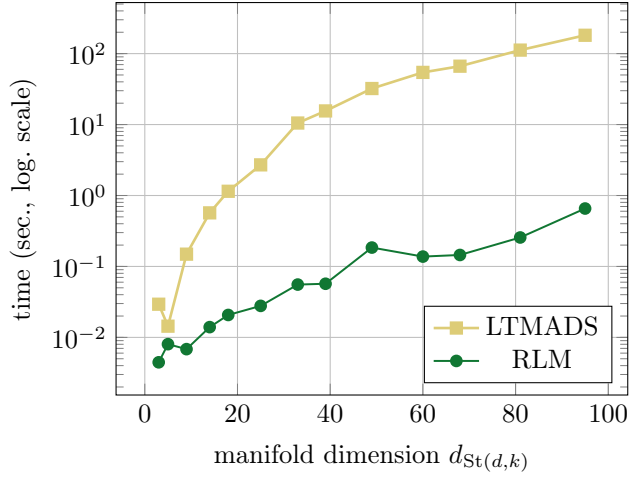
\begin{figure}[tbp]
    \centering
    \begin{tikzpicture}
        \begin{axis}[
            xlabel={manifold dimension $d_{\mathrm{St}(d,k)}$}, ylabel={time (sec., log.~scale)},
            legend pos=south east, grid = major,
            ymode = log,
            height=.33\textheight,
            width=.66\textwidth,
        ]
            \addplot[TolMutedSand, line width = 1pt, mark size = 2, mark=square*] table [x=dim, y=t2] {\Stdkdatatable};
            \addlegendentry{LTMADS}
            \addplot[TolMutedGreen, line width = 0.75pt, mark size = 2, mark=*] table [x=dim, y=t1] {\Stdkdatatable};
            \addlegendentry{RLM}
        \end{axis}
    \end{tikzpicture}
    \caption{Runtime comparison of our RLM (green) to the LTMADS algorithm (sand) on different Stiefel manifolds.}%
    \label{fig:Stdk-runtime}
\end{figure}

We again first compare the runtimes, cf.~Figure~\ref{fig:Stdk-runtime}.
We see a similar behaviour as for the previous experiment: RLM stays below $0.66$ seconds in runtime,
while the LTMADS increases up to about $180$ seconds,which is again a significant speed-up.

\begin{figure}[tbp]
    \centering
    \begin{tikzpicture}
    \begin{axis}[
        xlabel={manifold dimension $d_{\mathrm{St}(d,k)}$},
        ylabel={\# Iterations (log.~scale)},
        legend style={at={(0.97,0.5)},anchor=east},
        grid = major,
        ymode = log,
        height=.33\textheight,
        width=.66\textwidth,
    ]
        \addplot[TolMutedSand, , line width = 1pt, mark size = 2, mark=square*] table [x=dim, y=iter2] {\Stdkdatatable};
        \addlegendentry{LTMADS}
        \addplot[TolMutedGreen, , line width = 0.75pt, mark size = 2, mark=*] table [x=dim, y=iter1]{\Stdkdatatable};
        \addlegendentry{RLM}
        \end{axis}
    \end{tikzpicture}
    \caption{Cost value of the computed minimizers on $\mathrm{St}(d,k)$ of our RLM (green) to the LTMADS algorithm (sand).
    }%
    \label{fig:Stdk:Iterations}
\end{figure}
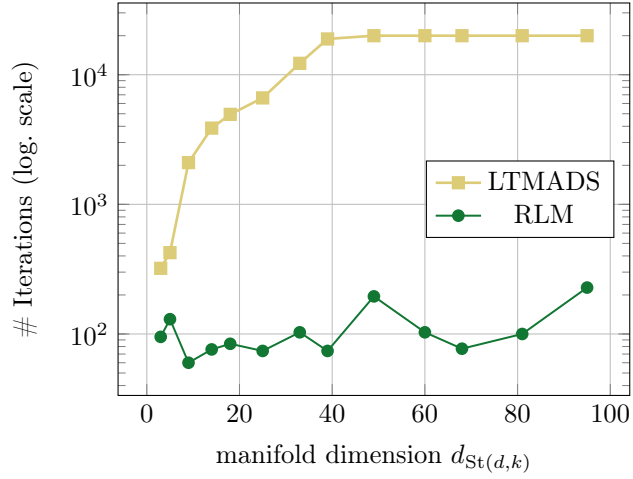
For the final cost value in this experiment, both RLM and LTMADS reach approximately the same value, where RLM is either numerically equal or ---in a few cases--- about $0.01\,\%$ better.
We additionally look at the number of iterations, each of the solver runs required, cf.~Figure~\ref{fig:Stdk:Iterations}.
Already for the smallest dimension $d=3$ LTMADS needs over 3 times more iterations, which afterwards increases very fast until we hit the maximal number of iterations from $d=11$ onwards. RLM on the other hand usually needs around 100 iterations, and in a few cases like for $d=11$ or $d=15$ it might increase up to about 200 iterations.

\subsection{Bundle adjustment with box constraints}

The bundle adjustment problem is a classical problem in computer vision and photogrammetry~\cite{Zach:2014}.
Here we extend the standard formulation with constraints on the camera parameters and point positions~\cite{GongMengSeibel:2015}.
All parts of the problem as defined below have been separately considered in the literature~\cite{GongMengSeibel:2015,Zach:2014}
,
but to the best of our knowledge, there is no prior work that supports robustification, optimization on manifolds and bounds constraints simultaneously.

The input data in bundle adjustment consists of $n_{\mathrm{obs}}$ observations $o \in (\mathbb{R}^2)^{n_{\mathrm{obs}}}$ describing 2D image coordinates of $n_{\mathrm{pts}}$ characteristic points in a scene, captured by $n_{\mathrm{cams}}$ cameras from different viewpoints.
The goal is to reconstruct 3D positions of points $p_{\mathrm{pts}} \in M_{\mathrm{pts}}$ in the scene, where
\begin{equation*}
    M_{\mathrm{pts}} = \left([l_{\mathrm{pts ,x}}, u_{\mathrm{pts, x}}] \times [l_{\mathrm{pts,y}}, u_{\mathrm{pts,y}}] \times [l_{\mathrm{pts,z}}, u_{\mathrm{pts,z}}]\right)^{n_{\mathrm{pts}}} \subset (\mathbb{R}^3)^{n_{\mathrm{pts}}}
\end{equation*}
and simultaneously estimate the camera positions $p_{\mathrm{t}} \in M_{\mathrm{t}}$, where
\begin{equation*}
    M_{\mathrm{t}} = \left([l_{\mathrm{t,x}}, u_{\mathrm{t,x}}] \times [l_{\mathrm{t,y}}, u_{\mathrm{t,y}}] \times [l_{\mathrm{t,z}}, u_{\mathrm{t,z}}]\right) ^{n_{\mathrm{cams}}} \subset (\mathbb{R}^{3})^{n_{\mathrm{cams}}},
\end{equation*}
orientations $p_{\mathrm{r}} \in \mathrm{SO}(3)^{n_{\mathrm{cams}}}$, and their intrinsic parameters: focal length $f \in [l_{\mathrm{f}}, u_{\mathrm{f}}]^{n_{\mathrm{cams}}}$ and radial distortion coefficients $k_1 \in [l_{\mathrm{k_1}}, u_{\mathrm{k_1}}]^{n_{\mathrm{cams}}}$, $k_2 \in [l_{\mathrm{k_2}}, u_{\mathrm{k_2}}]^{n_{\mathrm{cams}}}$, that best explain the observed 2D projections.
It is known that $i$th observation, $i\in \{1,\ldots,n_{\mathrm{obs}}\}$, was taken with camera $c(i)$ and corresponds to point $d(i)$, for some problem-dependent functions $c\colon \{1,\ldots,n_{\mathrm{obs}}\} \to \{1,\ldots,n_{\mathrm{cams}}\}$ and $d\colon \{1,\ldots,n_{\mathrm{obs}}\} \to \{1,\ldots,n_{\mathrm{pts}}\}$.
Point identification $d$ may contain mistakes, which can be considered as outliers in the data.
It is thus necessary to perform robust regression to obtain a reliable reconstruction of the scene and camera parameters.

The objective in the bundle adjustment problem reads
\begin{equation*}
    \argmin_{p_{\mathrm{t}}, p_{\mathrm{r}}, p_{\mathrm{pts}}, f, k_1, k_2} \sum_{i=1}^{n_{\mathrm{obs}}} \rho_i\bigl( \lVert F_{c(i)}(p_{\mathrm{t}, c(i)}, p_{\mathrm{r}, c(i)}, p_{\mathrm{pts}, d(i)}, f_{c(i)}, k_{1, c(i)}, k_{2, c(i)}) - o_i \rVert_2^2\bigr),
\end{equation*}
where $F_{c(i)}$ is the projection function that maps the 3D points and camera parameters, which reads
\begin{equation*}
    F_{c(i)}(p_{\mathrm{t}, c(i)}, p_{\mathrm{r}, c(i)}, p_{\mathrm{pts}, d(i)}, f_{c(i)}, k_{1, c(i)}, k_{2, c(i)})
    = f_{c(i)}\,\mathrm{rad}_{c(i)}\,
    \begin{pmatrix}
        x_{\mathrm{n},c(i)}\\
        y_{\mathrm{n},c(i)}
    \end{pmatrix},
\end{equation*}
where $x_{\mathrm{n},c(i)} = -\frac{x_{\mathrm{c},c(i)}}{z_{\mathrm{c},c(i)}}$, $y_{\mathrm{n},c(i)} = -\frac{y_{\mathrm{c},c(i)}}{z_{\mathrm{c},c(i)}}$ are the normalised image coordinates, $(x_{\mathrm{c},c(i)}, y_{\mathrm{c},c(i)}, z_{\mathrm{c},c(i)})$ are the coordinates of $p_{\mathrm{r},c(i)}\,p_{\mathrm{pts},d(i)} + p_{\mathrm{t},c(i)}$, $r_{c(i)}^2 = x_{\mathrm{n},c(i)}^2 + y_{\mathrm{n},c(i)}^2$ is the square of the radial distance, and $\mathrm{rad}_{c(i)} = 1 + k_{1,c(i)}\,r_{c(i)}^2 + k_{2,c(i)}\,r_{c(i)}^4$ is the radial distortion.

There are numerous software packages available for bundle adjustment, such as Ceres~\cite{CeresSolver:2023}, although they either don't support box constraints or optimization on manifolds.
The closest comparable implementation is the trust region reflective algorithm from SciPy~\cite{SciPy:2020}, which supports box constraints but is limited to Euclidean optimization.
Thus we resort to using the Euler angle parametrization of the rotation group $\mathrm{SO}(3)$ in SciPy for the camera orientation for the purpose of this experiment.
Moreover, since SciPy can only handle elementwise regularization instead of blockwise regularization like Ceres or our solver, in Python we put the robustifier directly in the residuals and don't use robustification on the solver side.
For our experiment we use Python 3.12 with SciPy 1.17.1.

For the numerical example we take the Ladybug dataset~\url{https://grail.cs.washington.edu/projects/bal/} with 49 images and subsample it to the first 20 cameras with corresponding point and iterations.
Initial rotations are set to the identity matrix, translations are set to 1, point positions to 0, focal lengths are set to 400 and radial distortion parameters to 0.
Furthermore, we use the Huber robustifier without scaling.
We constrain point position to the $[-1, 1]^3$ box,  focal lengths to the $[350, 450]$ interval and radial distortion parameters to the $[0, 0.1]$ interval.
Camera positions remain unconstrained.
We also set $\beta_{\mathrm{i}} = 8$, $\beta_{\mathrm{d}} = 0.2$, $\eta = 0.2$, $\eta_{\mathrm{l}} = 0.5$, $\mu_0 = 0.1$ and use strict scaling rule.

Figure~\ref{fig:BA-objective-history} shows the objective value history for both SciPy and our RLM implementation.
Although SciPy shows a faster decrease in the 10 iterations, our RLM implementation reaches a much lower objective value after about 100 iterations.
Both solvers take around \SI{0.09}{s} per iteration.

\pgfplotstableread[col sep=comma]{data/python_opt_history_subsampled.csv}\BAPythonHistory
\pgfplotstableread[col sep=comma]{data/julia_opt_history_subsampled.csv}\BAJuliaHistory
\begin{figure}[tbp]
    \centering
    \begin{tikzpicture}
        \begin{axis}[
            xlabel={iteration (log.~scale)},
            ylabel={objective value (log.~scale)},
            legend pos=north east,
            grid = major,
            xmode = log,
            xmin = 1,
            ymode = log,
            height=.33\textheight,
            width=.66\textwidth,
        ]
            \addplot[TolMutedSand, line width = 1.5pt, mark size = 0, mark=square*] table [x=iteration, y=objective] {\BAPythonHistory};
            \addlegendentry{SciPy TRF (Python)}
            \addplot[TolMutedGreen, line width = 1.5pt, mark size = 0, mark=*] table [x=iteration, y=objective] {\BAJuliaHistory};
            \addlegendentry{RLM (Julia)}
        \end{axis}
    \end{tikzpicture}
    \caption{Objective value history on the bundle-adjustment experiment for SciPy TRF (Python) and our RLM (Julia).}
    \label{fig:BA-objective-history}
\end{figure}
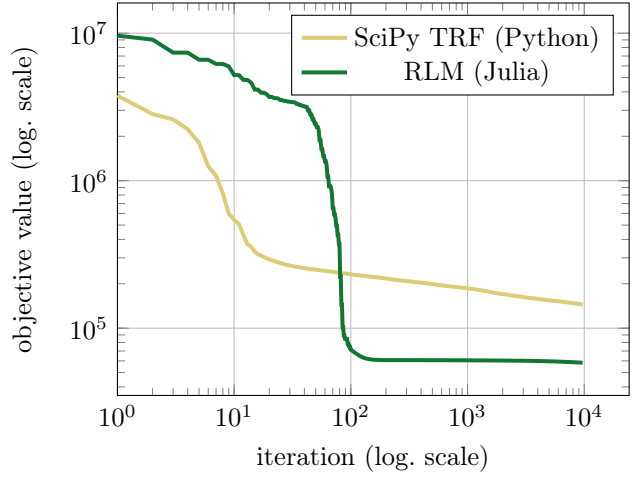

\section{Conclusion}
\label{section:conclusion}

In this paper generalized the robust Levenberg-Marquardt algorithm to Riemannian manifolds and operator-based Jacobians.
We thoroughly derived the algorithm, especially a new surrogate including second order terms of the robustifier.
This generalises the “Triggs correction” to Riemannian manifolds.
Even with the new surrogate and a slightly more general damping parameter update, the existing
convergence results also apply to this setting as well.
We further generalise the robust Levenberg–Marquardt algorithm to a Riemannian setting with Euclidean bound constraints.

An accompanying new release of \manoptjl provides a generic implementation of the introduced algorithm.
Using this open source code, we illustrate how a robust Levenberg-Marquardt qualitatively
improves geodesic regression under outliers. We further illustrate how the robust Procrustes problem
can be solved several orders of magnitude faster than with previous algorithm even when the measurements are incomplete.
Finally, the bundle adjustment problem was extended to include bound constraints. Here, our algorithm
outperforms existing software and finds a better minimizer.

\subsection*{Acknowledgements}

RB would like to thank MB and the AGH University for the kind hospitality and discussions during
a sabbatical stay in autumn 2025, where the main part of this paper stems from.
This research project supported by the program ,,Excellence initiative -- research university'' for the AGH University under the application IDUB 15636 (action D11).

\appendix

\printbibliography

\end{document}